\documentclass[12pt]{article}
\usepackage[amsmath]{e-jc}
\usepackage{amsfonts,relsize,physics,longtable,booktabs,array, algorithmic, mathrsfs,tikz,lmodern,diagbox,subfigure,caption}
\usepackage{algorithm2e}
\everymath{\displaystyle}
\usepackage{mathtools}
\usepackage{graphicx}
\newcommand{\Proba}{\mathbb{P}}
\newcommand{\R}{\mathbb{R}}
\newcommand{\N}{\mathbb{N}}
\newcommand{\D}{\text{d}}

\dateline{TBD}{TBD}{TBD}


\Copyright{The authors. Released under the CC BY-ND license (International 4.0).}

\title{An analytical framework for the Levine hats problem: \\ new strategies, bounds and generalizations.}

\author{Clément Bouquet*
\and
Salah Chikhi*
\and
Timoth\'e Charles
\and
Yanghao Zhou
\and
Eric Wang
}

\authortext{1}{
   \email{firstname.lastname@polytechnique.edu} \\
   *: Equal contribution
   }

\begin{document}

\maketitle
% ABSTRACT
% E-JC papers must include an abstract. The abstract should consist of a
% succinct statement of background followed by a listing of the principal
% new results that are to be found in the paper. The abstract should be
% informative, clear, and as complete as possible. Phrases like
% "we investigate..." or "we study..." should be kept to a minimum in
% favor of "we prove that..."  or "we show that...".  Do not include equation
% numbers, unexpanded citations (such as "[23]"), or any other references
% to things in the paper that are not defined in the abstract. The abstract
% may be distributed without the rest of the  paper so it must be entirely
% self-contained.  Try to include all words and phrases that someone
% might search for when looking for your paper.

\begin{abstract}
  We study the Levine hat problem, a cooperative puzzle introduced by Lionel Levine 
in 2010, in which $n \geq 2$ players must simultaneously identify a black hat on 
their own infinite stack, each seeing only their teammates' stacks. While the 
optimal winning probability $V_n$ remains unknown even for $n=2$, we make three 
key advances. First, we develop a geometric and integral framework representing 
strategies as Lebesgue-measurable functions, yielding a new integral expression 
for $V_n$ and a unified treatment of finite and infinite stacks. Second, we 
construct a recursive strategy $\mathscr{S}_5$ processing hats in blocks of five, 
which attains the conjectured optimal probability $7/20$ for two players. Although 
this bound was already achieved by the known strategy $\mathscr{S}_3$, the 
existence of $\mathscr{S}_5$ refutes the previously held expectation that recursive strategies 
with block size greater than three yield no improvement, and produces a strictly 
better geometric convergence rate for $V_{2,h}$ as well as a new lower bound for 
$V_2(p)$ which improves known results for $p < 0.312$. Building upon this, we improve the geometric convergence rate of $V_{2,h}$ up to the near-optimal $1/4^{1-\varepsilon}$ for any $\varepsilon > 0$.  Third, we introduce and completely solve a generalization of the problem where players are given uncountably infinite stacks of hats, showing that the optimal winning probability in this setting equals exactly $1/2$ for all $n \geq 2$. This new formulation allows to study the original combinatorial problem using tools from analytic optimization, and provides a natural framework for computing optimal responses to fixed strategies.
\end{abstract}

\newpage

\tableofcontents

\newpage

\section{Introduction}

\subsection{The usual definition of Levine's hat game}

In this section, we provide an overview of the topic covered in our work. We begin by presenting the problem in its original formulation. \\

In 2010, Lionel Levine conceived an elementary "hat" problem that would 
attract the attention of the mathematical research community. Although the 
problem may at first appear to be a purely recreational puzzle, it seems to 
have grown out of Levine's work on fast simulations of certain growth 
models~\cite{article 6}. It was later popularized by Tanya 
Khovanova~\cite{article 4} and became notorious for the difficulty of 
answering even its most basic questions despite its elementary formulation. 
Today it is commonly known as \textit{Levine's hat problem}. It is a cooperative game involving $n \geq 2$ players who each have a stack of $h$ hats on their heads. Each hat is either black or white. The colors of the hats are independently sampled according to a Bernoulli distribution with parameter $1/2$. Each player can then observe the stacks of the other players but cannot observe their own stack. During the game, no one is allowed to communicate, but the players can collectively devise a strategy beforehand. Each player then simultaneously selects a positive integer. The players win if and only if each player has selected the index of a black hat on their own head. 

The main questions of interest here are the following: with what probability can the players win? How does that probability evolve as $h$ goes to infinity? When $n$ goes to infinity?

\begin{figure}[h]
    \centering
    \subfigure[Losing game: the second player chooses index 1, whereas the first hat on his head is white.]{
        \includegraphics[width=0.45\textwidth]{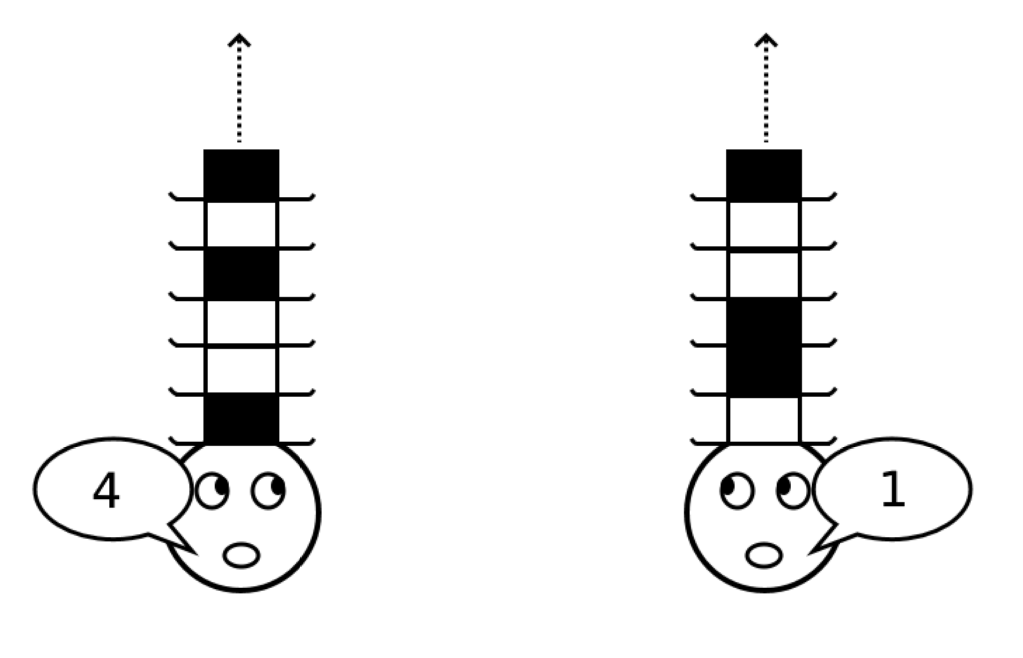}
    }
    \hspace{3pt}
    \subfigure[Winning game: Both players find the index of a black hat on their own heads.]{
        \raisebox{-7pt}{\includegraphics[width=0.4275\textwidth]{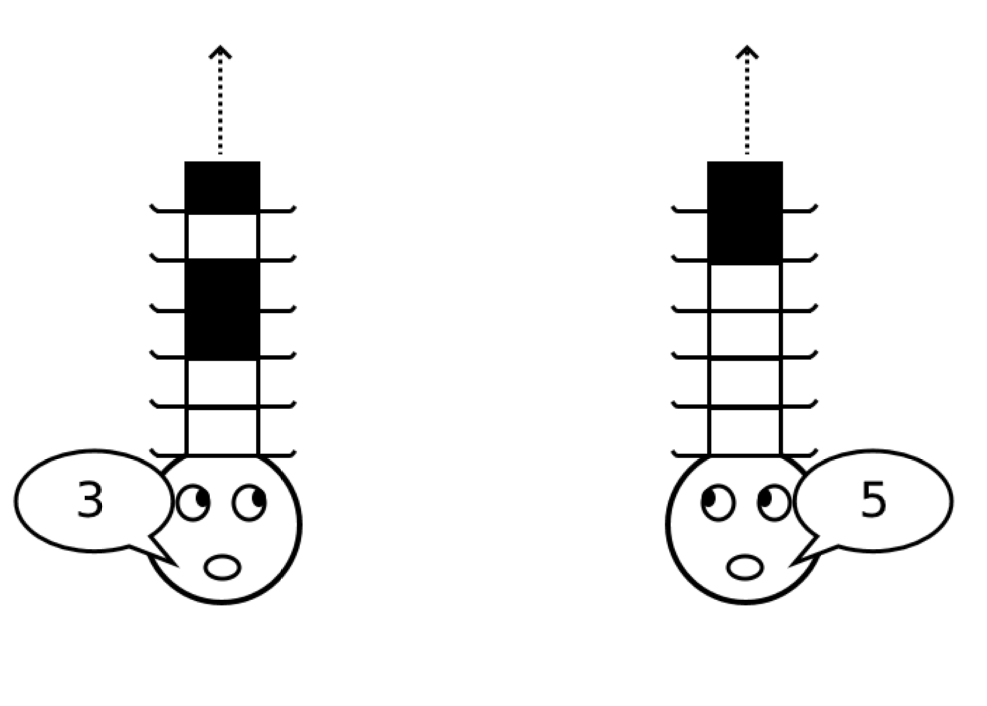}}
    }
    \caption{Illustration of Levine's hat problem ($n=2$)}
    \label{fig:deux_images}
\end{figure}

For small values of $h$, the problem can be entirely solved by brute force on a computer. However, the set of possible strategies becomes far too large very quickly to hope to achieve the same approach when $h$ becomes large. For that reason, the main case of interest is when $h$ tends to infinity. \\

We now give a few details to formalize what has just been stated. Each player $i$ receives a random stack denoted $U_{i} := \left( U_{i}^{(k)} \right)$, where the $U_{i}^{(k)}$ are independent random variables following a Bernoulli distribution with parameter $1/2$. Equivalently, one could also consider the uniform probability measure $\Proba$ on $\left( \{0,1\}^{h} \right)^{n}$. \\

We now define strategies on that set. As each player can observe the stacks of their teammates but not their own stack, an individual strategy is a function from $\left( \{0,1\}^{h} \right)^{n-1}$ to $\{1, \dots, h\}$. We will call this type of strategy an $h$-strategy (ignoring the dependence on $n$) with $h < \infty$. An $h$-strategy is thus a deterministic function giving the choice of a player depending on their teammates' stacks. We emphasize that there is no loss of generality in restricting to deterministic strategies, as these are just convex combinations of deterministic ones. Stated differently, for every stochastic strategy, there exists a deterministic strategy that is just as good. \\

From now on, we define $V_{n,h}$ to be the maximum probability of victory for Levine's hat problem where each player is given $h$ hats and the possible strategies are elements of $\mathcal{S}_{n-1,h}$, the set of $h$-strategies. We naturally take $\mathcal{S}_{0,h} = \{1, \dots h\}$: a strategy in the case where $n=1$ is constant. Equivalently: \begin{center}
    $V_{n,h} = \sup_{k_{1}, \dots , k_{n} \in \mathcal{S}_{n-1,h}} \Proba\left(\forall 1\leq i \leq n, U_{i}^{(k_{i}(U_{j}, j\neq i))} = 1\right)$
\end{center}

It follows that, for fixed $n$, the sequence $\left(V_{n,h}\right)_{h \geq 1}$ is increasing and trivially bounded above by $1/2$. Hence, it converges to a limit: \begin{center}
    $V_{n} := \lim_{h \to \infty} V_{n,h} \in \left] 0, \frac{1}{2} \right]$
\end{center}

This paper centers on the quantities $V_{n,h}$ and $V_{n}$. To be more precise, we will first present what makes this problem interesting and the conjectures associated with $\left(V_{n}\right)_{n \geq 1}$.

\subsection{Conjectures and previous works for Levine's hat problem}

At first glance, this problem might be counterintuitive. As the stacks are independent, the visible ones bring no extra information to a player's own hats. Hence, one may believe that all strategies achieve the same probability of victory since the hats are as likely to be black as to be white. More precisely, each hat is black with probability $1/2$. Hence, the probability for the team of players to win would be $1/2^{n}$ regardless of the strategy. However, it is possible to find strategies that have a probability of victory strictly larger than $1/2^{n}$. \\

We present a first example of such a strategy. Consider the case $n=2$ 
and $h=\infty$. Players A and B use the following strategy: each picks 
the index of the first black hat on their teammate's stack. Writing 
$m_A, m_B\geq 1$ for these indices, both players win if and only if 
$m_A = m_B$. At the index $\min(m_A,m_B)$, the three configurations 
$(1,0)$, $(0,1)$, $(1,1)$ are equally likely, and $(1,1)$ occurs 
if and only if $m_A=m_B$. Hence this ``first black hat'' (FBH) strategy 
achieves a winning probability of $1/3 > 1/4$, already beating the 
naive bound. We refer to Section~\ref{section algo} for a full 
analysis of this strategy within our framework.\\

This apparent paradox arises from the fact that the probability of winning associated with a family of $n$ individual strategies (one per player) involves correlations between the different strategies. The players can exploit these correlations by virtually reducing the configuration space, thereby increasing the proportion of outcomes that are collectively favorable to them. \\

Moreover, taking the limit in the previous example ($h \longrightarrow \infty$) makes it possible to disregard the boundary conditions specific to the case of finite stacks. However, this does not necessarily simplify our understanding of the problem, and several questions remain open concerning the quantity $V_n$. We now state the two main conjectures that motivate ongoing research on this problem and form the guiding thread of this paper.

\subsubsection{First conjecture}

The first question regarding Levine's game concerns the exact values of $V_n$. Surprisingly, this remains an open problem even in the case $n=2$. The best known strategy for the two-player version of Levine’s problem yields a winning probability of $7/20$. The first conjecture is therefore as follows.

\begin{conjecture}
    $V_2=7/20$
    \label{conj1}
\end{conjecture}

The existence of a strategy achieving a success probability of $7/20$ immediately establishes that $V_2 \geq 7/20$. To complete the proof, we must demonstrate the reverse inequality $V_2 \leq 7/20$. Several approaches have been developed for this purpose. In \cite{article 3}, the authors introduce the concept of "hints" provided to the players. Crucially, since players retain the freedom to disregard these hints, this modification can only improve their chances of success compared to the original game. Through exhaustive analysis of all possibilities in this enhanced framework -- where the hints effectively reduce the space of strategies requiring verification -- the authors establish that: \begin{center}
    $V_{2} \leq 81/224 \simeq 0.3616$
\end{center}

The more recent article \cite{article 5} provides a different approach. Indeed, using results from Fourier theory, it proves that: \begin{center}
    $V_{2} \leq 0.37193$
\end{center}

Although this method yields a weaker bound, it has the advantage of being derived entirely without computer assistance. \\

Those are the two main strategies developed by prior works for obtaining upper bounds. However, in the hope of obtaining more interesting results for this problem, similar studies have proposed variants of this puzzle. For instance, in section 2.2 of reference \cite{article 3}, Joe Buhler \textit{et al.} consider the two-player game with a biased coin. Stated differently, they study the same game where each hat has a probability $p \in [0,1]$ of being black. In particular, they provide lower bounds for the optimal probability of victory $V_{2}(p)$ which is now a function of $p$. In this new puzzle, Conjecture~\ref{conj1} can be rephrased as $V_{2}(1/2) = 0.35$. One of the main results they prove for this game is the following bound, which has yet to be improved :

\begin{enumerate}
    \item \hfill $\displaystyle U_{1}(p):= \frac{p(1 + p + p^2 + 3p^3 - 3p^4 + p^5)}{(1 + p)(2 - p)(1 + p^2)} \leq V_2(p)$ \hfill for $p \leq \frac{1}{2}$;
    
    \vspace{1em} % Adds a bit of breathing room between the two
    
    \item \hfill $\displaystyle U_{2}(p) := \frac{p(1 + 5p - 10p^2 + 10p^3 - 5p^4 + p^5)}{(2 - 2p + p^2)(1 + p)(2 - p)} \leq V_2(p)$ \hfill for $\frac{1}{2} \leq p$.
\end{enumerate}

Note that $U_{1}$ and $U_{2}$ represent the win probabilities of two distinct strategies within this modified game. The relative performance of these strategies depends strictly on the value of $p$: specifically, $U_{1}$ dominates $U_{2}$ on the interval $]0, 1/2[$, while $U_{2}$ becomes the superior strategy on $]1/2, 1[$. Note that $U_1(1/2) = U_2(1/2) = 7/20$, which is consistent with Conjecture \ref{conj1}.

\subsubsection{Second conjecture}

The second common question concerns the asymptotic behavior of the sequence $(V_n)_{n \geq 1}$ as $n$ tends to infinity. We can already state the following lemma.

\begin{lemma}
The sequence \( (V_n)_{n \geq 1} \) is non-increasing.
\end{lemma}

The key argument is the independence of the hats, which ensures that a given stack provides no information about the others. Thus, when ignoring the outcome of a fixed player in the team, one can also disregard that player's stack in the formulation of strategies without reducing the winning probability for the remaining players. It is noteworthy that the sequence $(V_{n})_{n \geq 1}$ is strictly decreasing, as demonstrated in Theorem 11 of section 4.2 in \cite{article 3} using combinatorial methods and in \cite{article 1} through graph-theoretic approaches.

It is clear that the winning probability in the game with a single player (\( n = 1 \)) is $1/2$. Since the sequence \( (V_n)_{n \geq 1} \) is non-increasing and bounded below by 0, it admits a limit in the interval $[0, 1/2]$. Other arguments refining these bounds can be found in the literature, but generally few results have been established. In particular, the value of this limit remains unknown. The best known lower bound was first proven by Peter Winkler in 2010. A proof of this result can be found in a reference \cite{article 3}:

\begin{proposition}
    There exists a constant \( C > 0 \) such that for all \( n \geq 2 \), 
\[
V_n \geq \frac{C}{\ln(n)}.
\]
\end{proposition}

This result is not sufficient to prove that $V_{n} \to 0$. However, if that convergence to $0$ holds, we can then study the speed of convergence which is at least logarithmic. In this regard, reference \cite{article 3} mentions an interesting fact: if we consider the version of the game where players are required to find the index of their $\textbf{first}$ black hat, then the optimal probability of victory tends to $0$. It is conjectured that $V_{n}$ has the same property in the initial version of the game. More precisely, the second main conjecture of this problem is:

\begin{conjecture}
    $V_{n} \to 0$
\end{conjecture}

Definitive proofs for these two conjectures have yet to be established. Nevertheless, the results derived thus far—bolstered by consistent computer-assisted simulations—point toward their veracity.

\subsection{Results and contributions}

In this subsection we summarize the main contributions of this paper. 
Faced with the long-standing impasse that has characterized this problem, we introduce a new angle of attack by bringing the tools of real analysis and Lebesgue integration to provide new insights on this purely combinatorial puzzle. Our work is organized around three central developments, each advancing the state of the art on this problem.

\medskip

\paragraph{A geometric and integral framework for strategies.}
We introduce a mathematical formalism in which player strategies are represented as measurable functions. 
By encoding infinite hat configurations as real numbers $x \in [0,1[$ via their binary expansion, we obtain a unified treatment of both finite and infinite stacks through Lebesgue integration. 
This representation admits a natural geometric interpretation: strategies correspond to black and white regions of the unit cube $[0,1]^n$, which can be visualized as a checkerboard. 
In particular, we show that $V_n$ admits an equivalent integral representation as the supremum of a Lebesgue integral over $[0,1]^n$, transforming the problem from a discrete combinatorial optimization into a variational one accessible to analytic methods.

\medskip

\paragraph{A new strategy and refined bounds.}
Building on this analytical framework, we construct in Theorem~\ref{s5} (Section~\ref{section algo})
a recursive strategy $\mathscr{S}_5$ of order $5$ achieving the conjectured optimal probability 
$7/20$ for two players and achieving better bounds for $V_{2,h}$ and $V_{2}(p)$. This refutes the conjecture of Buhler et al.\ \cite{article 3} that 
recursive strategies with block size $t > 3$ yield no improvement over those of order $3$. 
In Theorem~\ref{thm:Rt_exists}, we show that $\mathscr{S}_5$ belongs to an infinite family of strategies which we call $R_t$-type strategies. We show that 
for every odd $t \geq 3$, there exists an $R_t$-type strategy achieving $7/20$, which yields 
in Theorem~\ref{final_rate} a geometric convergence rate 
$V_{2,h} \geq 7/20 - C_\varepsilon / (4^{1-\varepsilon})^h$ for any $\varepsilon > 0$.  We also show that $R_t$-type strategies do not exist for even values of $t$.
Additionally, Theorem~\ref{u3(p)} establishes a new lower bound $U_3(p)$ for $V_2(p)$ 
which strictly improves all previously known bounds for $0 < p < 0.312$. All of these bounds originate from the discovery of $\mathscr{S}_{5}$.

\medskip

\paragraph{Continuous generalization and imaginary strategies.}
We introduce the \emph{continuous Levine game}, in which players choose 
hat indices from $\mathbb{R}_+$ rather than $\mathbb{N}^*$, giving rise 
to a class of strategies which we call \emph{imaginary strategies}. We 
prove that the optimal winning probability satisfies $W_n = 1/2$ for all 
$n \geq 2$ (Theorem~\ref{thm W_n >= V_n}). Unlike previous variants which 
relax the rules or information structure of the game, our generalization 
relaxes the \emph{data structure} itself — the discrete binary encoding of 
hat stacks — and the result identifies this binary constraint as the 
essential obstruction keeping $V_n$ below $1/2$: once lifted, perfect 
correlation becomes achievable and the problem reduces to a tractable 
analytic one. Within the mixed framework $2^{\widehat{\mathcal{S}}_1} 
\times \widehat{\mathcal{M}}_1$, we establish an Optimal Response Theorem 
(Theorem~\ref{optimresp}) and show that the first-black-hat strategy is 
its own best response, providing first evidence for 
Conjecture~\ref{new conjecture} and a potential path toward a sharp bound 
on $V_2$.
\section{Framework unifying finite and infinite stacks}
\label{section visualisation des stratégies utilisées par les joueurs}

We fix $n=2$. Recall that an $h$-strategy is a map $k:\{0,1\}^h \to \{1,\ldots,h\}$.
We index the $2^h$ possible stacks lexicographically: $a_1 < \cdots < a_{2^h}$
in $\{0,1\}^h$.\\

\begin{figure}[h!]
\centering
\begin{tabular}{cc}
\includegraphics[width=0.35\linewidth]{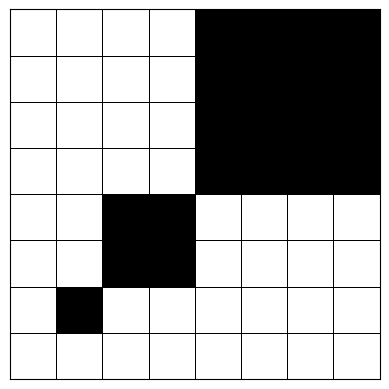} &
\includegraphics[width=0.35\linewidth]{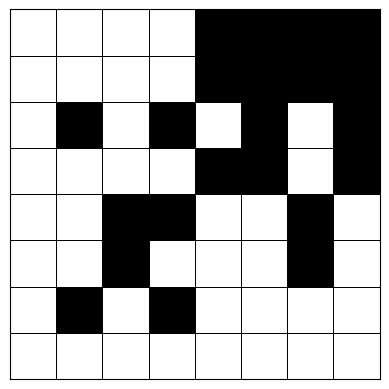} \\
$\delta^{\mathscr{S}_{\mathrm{FBH},3},\mathscr{S}_{\mathrm{FBH},3}}$,
winning prob.\ $21/64$ &
$\delta^{\mathscr{S}_{3,3},\mathscr{S}_{3,3}}$,
winning prob.\ $22/64$
\end{tabular}
\caption{Joint outcome matrices for $h=3$, $n=2$.
Left: first-black-hat strategy $\mathscr{S}_{\mathrm{FBH},3}$.
Right: strategy $\mathscr{S}_{3,3}$, optimal for $h=3$
as verified by brute force.
The fractal structure visible in the $\infty$-strategy limits
(Figures~\ref{fig:strategies-comparison1} and~\ref{fig:strategies-comparison})
is already suggested here.}
\label{canonicalhstrategies}

\end{figure}

Given strategies $k_1, k_2$ for players $A$ and $B$, define the
\emph{joint outcome matrix} $\delta^{k_1,k_2}:\{1,\ldots,2^h\}^2\to\{0,1\}$ by
\[
\delta^{k_1,k_2}(i,j) := a_i^{(k_1(a_j))}\cdot a_j^{(k_2(a_i))},
\]
where $a_i^{(k)}$ denotes the $k$-th bit of $a_i$. The entry
$\delta^{k_1,k_2}(i,j)=1$ if and only if both players simultaneously
select a black hat in configuration $(i,j)$. Since all $4^h$ configurations
are equally likely,
\[
\Proba\!\left(U^{(k_1(V))}=V^{(k_2(U))}=1\right)
= \frac{1}{4^h}\sum_{i,j}\delta^{k_1,k_2}(i,j)
=: \mathcal{A}_{k_1,k_2},
\]
which equals the fraction of black tiles when $\delta^{k_1,k_2}$ is rendered
as a $2^h\times 2^h$ checkerboard (column $i$ left to right,
row $j$ bottom to top). We use this visualization extensively throughout
the paper. Figure~\ref{canonicalhstrategies} gives the representations of two canonical strategies: the FBH strategy and $\mathscr{S}_{3,3}$. $\mathscr{S}_{3,3}$ was first introduced in \cite{article 3} and serves as their starting point for a strategy achieving win probability of $7/20$. \\

Building upon this finite case, we now rigorously define Levine's hat problem with infinite stacks,
using the Lebesgue measure. This will allow us to unify these two scenarios. The key idea is to identify a stack (a sequence of $0$s and $1$s) with a real number $x\in[0,1[$ via
its binary expansion. This yields a new integral expression for $V_n$,
a natural extension of the visualization tool, and a framework
in which the best known strategies arise as measurable functions.

\subsection{From infinite stacks to Lebesgue measure}

An infinite random stack is a sequence $(U^{(k)})_{k\geq 1}$ of independent
Bernoulli$(1/2)$ variables. The natural probability space is
$\Omega:=\{0,1\}^{\mathbb{N}^*}$, equipped with the $\sigma$-algebra
$\mathcal{A}$ generated by cylinders
$C(\varepsilon):=\{u\in\Omega\mid u^{(i)}=\varepsilon^{(i)},\,1\leq i\leq m\}$
and the unique measure $\Proba$ satisfying $\Proba(C(\varepsilon))=2^{-m}$. \\

Stacks are identified with real numbers via the bijection
\[
\varphi:\Omega\setminus\mathcal{N}\longrightarrow[0,1[,\quad
u\longmapsto\sum_{i=1}^{\infty}\frac{u^{(i)}}{2^i},
\]
where $\mathcal{N}$ is the countable set of eventually-$1$ sequences. In what follows, we let $\Lambda$ be the Lebesgue measure on $\mathbb{R}$.

\begin{theorem}
\label{unique measure theorem}
The map $\Proba:=\Lambda\circ\varphi$ is the unique probability measure on
$\mathcal{A}\setminus\mathcal{N}$ satisfying $\Proba(C(\varepsilon))=2^{-m}$.
Furthermore, $\varphi$ is a measure-preserving bijection from
$(\Omega\setminus\mathcal{N},\mathcal{A}\setminus\mathcal{N},\Proba)$
to $([0,1[,\mathcal{B},\Lambda)$.
\end{theorem}

\subsection{$\infty$-strategies and integral formulation}

We introduce the following notation throughout:
\begin{itemize}
\item $\bar\Omega := \Omega\setminus\mathcal{N}$,
  and $\bar\varphi_m:\bar\Omega^m\to[0,1[^m$ denotes
  the coordinatewise application of $\varphi$.
\item $(x^{(i)})_{i\geq 1}$ denotes the binary expansion of $x\in[0,1[$.
\item $\varepsilon(x):=\lfloor x\rfloor-2\lfloor x/2\rfloor$,
  so that $x^{(i)}=\varepsilon(2^i x)$ for all $x\in[0,1[$ and $i\geq 1$.
\item $x_{-i}:=(x_1,\ldots,x_{i-1},x_{i+1},\ldots,x_n)$.
\item $\widehat{\mathcal{S}}_m$ denotes the set of Lebesgue-measurable
  functions from $[0,1[^m$ to $\mathbb{N}^*$.
\end{itemize}

We call any element of $\widehat{\mathcal{S}}_{n-1}$ an
\emph{$\infty$-strategy} for $n$ players. 
The identification of stacks with real numbers via $\varphi$ allows us to 
transfer measurability questions between the two spaces $\bar\Omega^m$ and 
$[0,1[^m$. The following lemma makes this transfer precise, and will be used 
repeatedly in what follows.

\begin{lemma}
\label{f o phi-1 measurable}
Let $m\geq 1$ and let $f:\bar\Omega^m\to\mathbb{R}$.
Then $f\circ\bar\varphi_m^{-1}$ is measurable if and only if $f$
is measurable. In particular, the map
$k\mapsto k\circ\bar\varphi_m^{-1}$ is a bijection from
$\mathcal{S}_m$ to $\widehat{\mathcal{S}}_m$, with inverse
$\hat k\mapsto\hat k\circ\bar\varphi_m$.
\end{lemma}

\begin{proof}
Both directions follow from the fact that $\bar\varphi_m$ and
$\bar\varphi_m^{-1}$ are measurable (Theorem~\ref{unique measure theorem}),
and that compositions of measurable functions are measurable.
The bijection statement follows immediately.
\end{proof}

Before writing down the integral expression for the winning probability, 
we need to verify that the event of winning is measurable in the first place 
when players use $\infty$-strategies. This is not immediate: the event 
$\left(\forall i,\, U_i^{(\hat k_i(U_{-i}))}=1\right)$ involves evaluating 
the $\hat k_i(U_{-i})$-th bit of $U_i$, which is a complex composition of 
measurable functions. The key is to rewrite this event entirely in terms of 
the function $\varepsilon$ and the bijection $\bar\varphi_n$, after which 
Lemma~\ref{f o phi-1 measurable} allows us to conclude.

\begin{lemma}
\label{lem:measurability}
Let $\hat k_1,\ldots,\hat k_n\in\widehat{\mathcal{S}}_{n-1}$.
Then the event $\left(\forall i,\, U_i^{(\hat k_i(U_{-i}))}=1\right)$
is measurable.
\end{lemma}

\begin{proof}
We compute:
\begin{align*}
\left(\forall i,\, U_i^{(k_i(U_{-i}))}=1\right)\cap\bar\Omega^n
&= \bar\varphi_n^{-1}\!\left\{(x_1,\ldots,x_n)\in[0,1[^n \;\middle|\;
   \forall i,\,\varepsilon\!\left(2^{\hat k_i(x_{-i})}x_i\right)=1\right\}.
\end{align*}
Since $\varepsilon$ and each $\hat k_i$ are measurable, the set inside
$\bar\varphi_n^{-1}$ is measurable, and $\bar\varphi_n^{-1}$ is measurable
by Lemma~\ref{f o phi-1 measurable}. Adding the negligible set $\mathcal{N}^n$
yields measurability of the full event.
\end{proof}

With measurability established, we can now compute the winning probability 
explicitly. The measure-preserving property of $\bar\varphi_n$ 
(Theorem~\ref{unique measure theorem}) allows us to transport the probability 
measure $\Proba$ to the Lebesgue measure on $[0,1[^n$, turning the winning 
probability into a Lebesgue integral. The identity $x^{(k)} = \varepsilon(2^k x)$ 
then yields the following expression.

\begin{lemma}
\label{lem:integral-rep}
For any $\hat k_1,\ldots,\hat k_n\in\widehat{\mathcal{S}}_{n-1}$,
\[
\Proba\!\left(\forall i,\, U_i^{(\hat k_i(U_{-i}))}=1\right)
= \int_{[0,1]^n}
  \prod_{i=1}^n\varepsilon\!\left(2^{\hat k_i(x_{-i})}x_i\right)
  \mathrm{d}\Lambda^n =: I_n.
\]
\end{lemma}

\begin{proof}
Using the measure-preserving bijection $\bar\varphi_n$ and the identity
$x^{(k)}=\varepsilon(2^k x)$:
\begin{align*}
\Proba\!\left(\forall i,\, U_i^{(\hat k_i(U_{-i}))}=1\right)
&= \Lambda^n\!\left(\bar\varphi_n\!\left(
   \left\{(u_j)\in\bar\Omega^n\mid\forall i,\,
   u_i^{(\hat k_i(u_{-i}))}=1\right\}\right)\right) \\
&= \int_{[0,1]^n}
   \prod_{i=1}^n x_i^{(\hat k_i(x_{-i}))}\,\mathrm{d}\Lambda^n
= \int_{[0,1]^n}
   \prod_{i=1}^n\varepsilon\!\left(2^{\hat k_i(x_{-i})}x_i\right)
   \mathrm{d}\Lambda^n. \qedhere
\end{align*}
\end{proof}

Lemmas~\ref{f o phi-1 measurable}--\ref{lem:integral-rep} together show 
that the supremum over $\infty$-strategies of the winning probability 
can be written as a Lebesgue integral. It remains to show that this 
supremum equals $V_n$, which was defined as $\lim_{h\to\infty} V_{n,h}$.
The inequality $I_n\geq V_n$ is straightforward: any $h$-strategy can be 
viewed as an $\infty$-strategy via $\varphi$, so the supremum over 
$\infty$-strategies is at least $V_{n,h}$ for every $h$.

The reverse inequality $I_n \leq V_n$ is more delicate. It requires showing 
that any $\infty$-strategy $\hat k$ can be approximated arbitrarily well 
by $h$-strategies, in the sense that the corresponding winning probabilities 
converge. The natural approximation is to replace $\hat k$ by the 
$h$-strategy $\hat k^{(h)}$ that agrees with $\hat k$ on each dyadic 
interval $I_{h,r}$ where $\hat k$ is ``almost constant''. The following 
lemma shows that this approximation converges pointwise almost everywhere, 
which by dominated convergence is enough to conclude.

\begin{lemma}
\label{lem:approx}
For any $\hat k\in\widehat{\mathcal{S}}_1$, the sequence of $h$-strategies
\[
\hat k^{(h)} := \sum_{i=1}^{h} i \sum_{r=0}^{2^h-1}
\delta\!\left(I_{h,r}\subset\hat k^{-1}(\{i\})\right)\mathbf{1}_{I_{h,r}}
\]
satisfies $\hat k^{(h)}\to\hat k$ pointwise almost everywhere as $h\to\infty$.
\end{lemma}

\begin{proof}
Since $\hat k$ is measurable, for each $i\geq 1$ the set
$\hat k^{-1}(\{i\})$ is measurable, and Lebesgue's measure regularity gives
\[
\bigcup_{\substack{0\leq r\leq 2^h-1\\
I_{h,r}\subset\hat k^{-1}(\{i\})}}I_{h,r}
\;\xrightarrow{\Lambda}_{h\to\infty}\;\hat k^{-1}(\{i\}).
\]
The family of such unions is increasing in $h$: for $h'>h$, each $I_{h,r}$
is a disjoint union of intervals $I_{h',r'}$, so
\[
\bigcap_{h'\geq h}
\bigcup_{\substack{0\leq r'\leq 2^{h'}-1\\
I_{h',r'}\subset\hat k^{-1}(\{i\})}}I_{h',r'}
= \bigcup_{\substack{0\leq r\leq 2^h-1\\
I_{h,r}\subset\hat k^{-1}(\{i\})}}I_{h,r}.
\]
Since $\hat k$ is defined a.e.,
$[0,1]\overset{\Lambda}{=}\bigcup_{i\geq 1}\bigcup_{h\geq 1}
\bigcap_{h'\geq h}\bigcup_{I_{h',r'}\subset\hat k^{-1}(\{i\})}I_{h',r'}$.
Hence for a.e.\ $x$, there exist $i\geq 1$ and $H\geq 1$ such that
$\hat k^{(h)}(x)=i=\hat k(x)$ for all $h\geq H$.
\end{proof}
We can now state and prove the main theorem of this section.

\begin{theorem}
\label{theoreme de reformulation complète}
For all $n\geq 1$,
\[
V_n
= \sup_{\hat k_1,\ldots,\hat k_n\in\widehat{\mathcal{S}}_{n-1}}
\int_{[0,1]^n}
\prod_{i=1}^n\varepsilon\!\left(2^{\hat k_i(x_{-i})}x_i\right)
\mathrm{d}\Lambda^n.
\]
\end{theorem}

\begin{proof}
Denote the right-hand side by $I_n$.
By Lemma~\ref{lem:integral-rep}, taking the supremum over
$\widehat{\mathcal{S}}_{n-1}$ and using the bijection of
Lemma~\ref{f o phi-1 measurable}, we have
\[
I_n = \sup_{k_1,\ldots,k_n\in\mathcal{S}_{n-1}}
\Proba\!\left(\forall i,\, U_i^{(k_i(U_{-i}))}=1\right).
\]

\noindent$(I_n\geq V_n)$: For each $h\geq 1$, let $\widehat{\mathcal{S}}_{m,h}$
denote the functions in $\widehat{\mathcal{S}}_m$ that are constant on
every dyadic interval $I_{h,r}$; these correspond exactly to $h$-strategies
via $\varphi$. Since $\widehat{\mathcal{S}}_{m,h}\subset\widehat{\mathcal{S}}_m$,
we have $I_n\geq V_{n,h}$ for all $h$, hence $I_n\geq V_n$.\\

\noindent$(I_n\leq V_n)$: We detail $n=2$; the case $n\geq 3$ follows
coordinate-wise. Fix $\hat k_1,\hat k_2\in\widehat{\mathcal{S}}_1$ and
let $\hat k_1^{(h)},\hat k_2^{(h)}$ be the approximating sequences
from Lemma~\ref{lem:approx}. By that lemma,
$(x,y)\mapsto\varepsilon(2^{\hat k_j^{(h)}(y)}x)$
converges pointwise a.e.\ to
$(x,y)\mapsto\varepsilon(2^{\hat k_j(y)}x)$
for $j=1,2$.
Since each factor is dominated by $1$, the dominated convergence
theorem gives:
\[
\int_{[0,1]^2}\varepsilon\!\left(2^{\hat k_1^{(h)}(y)}x\right)
\varepsilon\!\left(2^{\hat k_2^{(h)}(x)}y\right)\mathrm{d}x\,\mathrm{d}y
\;\longrightarrow\;
\int_{[0,1]^2}\varepsilon\!\left(2^{\hat k_1(y)}x\right)
\varepsilon\!\left(2^{\hat k_2(x)}y\right)\mathrm{d}x\,\mathrm{d}y.
\]
Since $(\hat k_1^{(h)},\hat k_2^{(h)})\in\widehat{\mathcal{S}}_{1,h}^2$,
the left-hand side is at most $V_{n,h}$ for each $h$.
Taking $h\to\infty$ gives
$\int_{[0,1]^2}\varepsilon(2^{\hat k_1(y)}x)
\varepsilon(2^{\hat k_2(x)}y)\,\mathrm{d}x\,\mathrm{d}y\leq V_n$.
Since $\hat k_1,\hat k_2$ were arbitrary, $I_n\leq V_n$.
\end{proof}
\subsection{Visualization of $\infty$-strategies}

For $n=2$ and strategies $k_1,k_2\in\widehat{\mathcal{S}}_1$,
the winning probability is
\[
P_{k_1,k_2}=\int_{[0,1]^2}
\varepsilon\!\left(2^{k_1(y)}x\right)\varepsilon\!\left(2^{k_2(x)}y\right)
\mathrm{d}x\,\mathrm{d}y,
\]
which equals the Lebesgue measure of the black region in $[0,1]^2$.
The following theorem confirms that this continuous visualization is
consistent with the finite one of Section~\ref{section visualisation
des stratégies utilisées par les joueurs}.

\begin{theorem}
\label{prolongement visualisation}
The $\infty$-strategy representation naturally extends the $h$-strategy
checkerboard visualization.
\end{theorem}

\begin{proof}
If $\hat k$ is constant equal to $K\leq h$ on $I_{j,h}=[j/2^h,(j+1)/2^h[$,
then for all $(x,y)\in I_{i,h}\times I_{j,h}$,
$\varepsilon(2^{\hat k(y)}x)=\varepsilon(2^K x)=x^{(K)}=a_i^{(K)}$.
The partition $(I_{i,h}\times I_{j,h})$ of $[0,1[^2$ into $4^h$ squares
then coincides exactly with the finite checkerboard, with matching indices.
\end{proof}

\subsection{Recursive strategies}

The best known strategies for Levine's problem share a common recursive
structure, which we now formalize.

\begin{definition}
\label{def:rec}
A strategy $K$ is \emph{recursive of order $h$} if it is defined by the
following algorithm, given a non-empty set
$\mathcal{C}\subsetneq\{0,1\}^h$ and an $h$-strategy $k_0$:
inspect the observed stack in consecutive blocks of $h$ hats;
skip each block not in $\mathcal{C}$; upon reaching the first block
$u\in\mathcal{C}$ after $m\geq 0$ skipped blocks, return $mh+k_0(u)$.
We call $k_0$ the \emph{associated $h$-strategy} of $K$.
\end{definition}

\begin{lemma}
Recursive strategies are $\infty$-strategies that terminate almost surely.
\end{lemma}

\begin{proof}
We prove the two claims separately.

\medskip
\noindent\textit{Termination almost surely.}
Fix a recursive strategy $K$ of order $h$ with stopping set 
$\mathcal{C}\subsetneq\{0,1\}^h$. For $m\geq 1$, let $A_m$ denote the 
event that the $m$-th block of $h$ hats belongs to $\mathcal{C}$. Since 
the hats are independent and $\mathcal{C}$ is non-empty, we have 
$\Proba(A_m) = |\mathcal{C}|/2^h > 0$ for all $m$, and the events 
$(A_m)_{m\geq 1}$ are independent. Since 
$\sum_{m=1}^\infty \Proba(A_m) = +\infty$,
the Borel--Cantelli lemma implies that $A_m$ occurs for infinitely many 
$m$ almost surely. In particular, the algorithm terminates almost surely. 

\medskip
\noindent\textit{Measurability.}
We show that $K\in\widehat{\mathcal{S}}_1$, i.e., that $K$ is a 
measurable function from $(\bar\Omega,\bar{\mathcal{A}},\Proba)$ to 
$\mathbb{N}^*$. It suffices to show that $K^{-1}(\{i\})$ is measurable 
for every $i\in\mathbb{N}^*$.

Fix $i\in\mathbb{N}^*$ and write $m = \lfloor(i-1)/h\rfloor$ and 
$j = i - mh\in\{1,\ldots,h\}$. The event $K^{-1}(\{i\})$ consists 
of all stacks $u\in\bar\Omega$ such that:
\begin{enumerate}
    \item the first $m$ blocks $u_1,\ldots,u_m$ of $h$ hats each lie 
    outside $\mathcal{C}$, and
    \item the $(m+1)$-th block $u_{m+1}$ lies in $\mathcal{C}$ and 
    satisfies $k_0(u_{m+1}) = j$.
\end{enumerate}
Each of these conditions concerns a fixed finite set of coordinates of 
$u$ and thus defines a cylinder of $\Omega$. More precisely:
\begin{itemize}
    \item The condition ``block $\ell$ lies outside $\mathcal{C}$'' is 
    
    $\left\{u\in\Omega\mid (u^{((\ell-1)h+1)},\ldots,u^{(\ell h)})
    \notin\mathcal{C}\right\}$,
    which is a finite union of cylinders since $\{0,1\}^h\setminus\mathcal{C}$ 
    is finite.
    \item The condition ``block $m+1$ lies in $\mathcal{C}$ and 
    $k_0(u_{m+1})=j$'' is similarly a finite union of cylinders, 
    since $\mathcal{C}$ and $k_0^{-1}(\{j\})\cap\mathcal{C}$ are 
    both finite subsets of $\{0,1\}^h$.
\end{itemize}
The set $K^{-1}(\{i\})$ is therefore a finite intersection of finite 
unions of cylinders, hence a finite union of cylinders. Since every 
cylinder belongs to $\mathcal{A}$, we conclude that $K^{-1}(\{i\})$ 
is measurable. By Lemma~\ref{f o phi-1 measurable}, the corresponding 
$\infty$-strategy $\hat K = K\circ\bar\varphi^{-1}$ is an element of 
$\widehat{\mathcal{S}}_1$.
\end{proof}

We emphasize that all currently known strategies achieving $7/20$ are recursive (or are combinations of recursive strategies), and the two canonical examples are as follows.

\begin{figure}[h!]
\centering
\begin{tabular}{cc}
\includegraphics[width=0.33\linewidth]{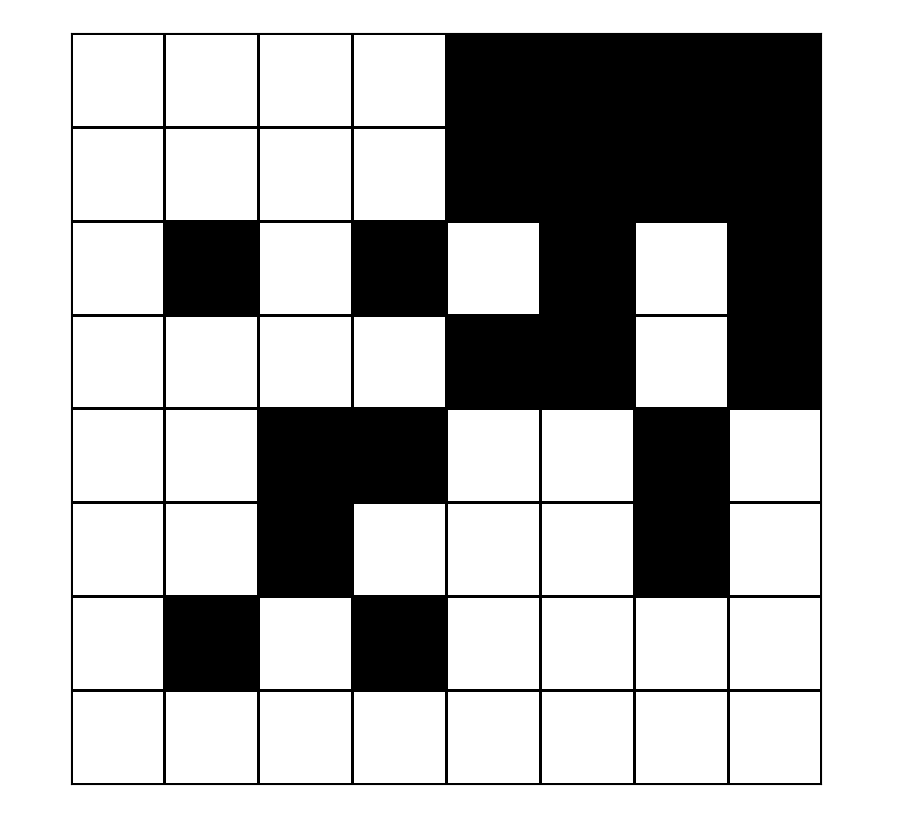} &
\includegraphics[width=0.3\linewidth]{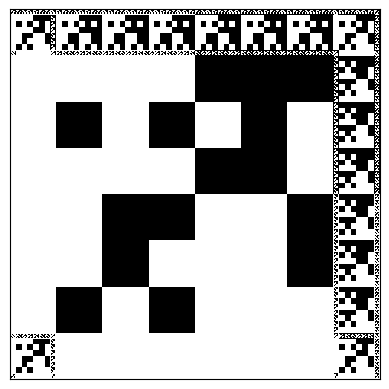} \\
3-strategy $\mathscr{S}_{3,3}$ &
Strategy $\mathscr{S}_3$ (winning prob.\ $7/20$)
\end{tabular}
\caption{Left: joint outcome matrix of $\mathscr{S}_{3,3}$ for $h=3$
(winning probability $22/64$).
Right: joint outcome matrix of the recursive $\infty$-strategy $\mathscr{S}_3$
(winning probability $7/20$).
The fractal structure arises from the recursive skipping of monochromatic
triplets (Definition~\ref{def:rec}).}
\label{fig:strategies-comparison1}
\end{figure}

\paragraph{First black hat strategy $\mathscr{S}_{\mathrm{FBH}}$.}
Taking $\mathcal{C}=\{(1)\}$ and $k_0=1$ is equivalent to the strategy where
each player names the index of the first black hat on their teammate's head.
The winning probability is $\sum_{m=1}^\infty 1/4^m = 1/3$.
Indeed, geometrically, the joint outcome matrix consists of a sequence of squares
along the diagonal with areas $1/4^m$ (see Figure~\ref{fig:strategies-comparison1}).

\paragraph{Strategy $\mathscr{S}_3$.}
Take $\mathcal{C}$ to be the set of non-monochromatic triplets and
$k_0=\mathscr{S}_{3,3}$, where
\[
(k_1(a_j))_{1\leq j\leq 8}=(1,3,2,2,1,3,1,1),\quad
(k_2(a_j))_{1\leq j\leq 8}=(1,3,2,3,1,1,2,1).
\]
In other words, each player advances over monochromatic triplets until reaching a
non-monochromatic one, then applies $\mathscr{S}_{3,3}$. This is one of the core strategies developed in \cite{article 3} which achieves a winning probability of $7/20$. \\

Among recursive strategies, a subclass of particular importance will
emerge in Section~\ref{section algo}: that of standardized strategies
(Definition~\ref{def:standardized}), for which the response to any
monochromatic block is independent of its color.
The winning probability of $\mathscr{S}_3$ equals $7/20$,
and its fractal visualization is shown in
Figure~\ref{fig:strategies-comparison1}.

\section{The strategy $\mathscr{S}_5$ and its consequences}
\label{section algo}

In this section, we present important results obtained algorithmically. These algorithmic results are the starting point for the discovery of new bounds.

\subsection{Discovery of efficient $h$-strategies}

We restrict ourselves to the case where stacks are finite and contain $h$ hats. Multiple algorithmic methods are possible to approach this problem. A direct brute-force approach is sufficient to find optimal strategies for values of $h \leq 3$. Beyond this threshold, the set of possible $h$-strategies becomes far too large to be explored entirely in a reasonable amount of time. For this reason, we decided to use a hill-climbing algorithm along with some heuristics with the goal of finding efficient $h$-strategies for larger values of $h$. \\

The hill-climbing algorithm proceeds as follows: starting from a random h-strategy, at each step we randomly select a stack configuration $a_j$ and update $k(a_j)$ randomly. The algorithm terminates when no single-coordinate update improves the objective. We also used simulated annealing which provided small adjustments that improved the bounds given by hill-climbing. The code used to generate this table is available \href{https://github.com/RayaneChikhi/An_analytical-framework_for-the_Levine-hats_problem/tree/master}{here}. \\

The hill-climbing algorithm allowed us to establish the following table of values. For each value of $h$, we provide a lower bound for $V_{2,h}$. This lower bound is obtained by finding an $h$-strategy that achieves such a probability of victory. However, it is likely that these bounds are sharp for $4 \leq h \leq 6$. Indeed, extensive attempts did not yield any better results.

    \begin{longtable}{|c|c|p{8cm}|}
    \hline
    \textbf{h} & \textbf{$V_{2,h} \geq \hdots$}   \\
    \hline
    4 & 0.34765625 
    \\
    \hline
    5 & 0.349609375  \\
    \hline
    6 & 0.349853515625  \\
    \hline
    7 & 0.34991455078125  \\
    \hline

    8 & 0.3499603271484375   \\
    \hline
    9 & 0.3499794006347656  \\
    \hline 

    10 & 0.34998035430908203  \\
    \hline
    \end{longtable}

This table of values clearly supports Lionel Levine's conjecture. Note that we lose significant potential by restricting ourselves to a finite number of hats. However, as seen in the previous section, we can construct $\infty$-strategies starting from an $h$-strategy. In particular, we can use the previous algorithmic results to create recursive $\infty$-strategies. Surprisingly, these will enable us to prove a novel lower bound estimations for the sequence $(V_{2,h})_{h \geq 1}$.

\subsection{A new recursive strategy of order 5}

In this subsection, we use the previous algorithmic results to construct a recursive $\infty$-strategy that achieves a probability of victory of $0.35$. As we will see, the mere existence of such a strategy is enough to establish two novel results on the problem.

\subsubsection{Definition of the strategy $\mathscr{S}_{5}$}

\begin{theorem}
\label{s5}
For $n=2$, there exists a symmetrical recursive $\infty$-strategy of order 5 that achieves a probability of victory equal to $7/20 = 0.35$.
\end{theorem}

\begin{proof}
We provide an example of such a strategy (which we construct using a $5$-strategy $\mathscr{S}_{5,5}$). Using the notations from the definition of recursive strategies, each player uses the recursive strategy defined by $\mathcal{C}$ and $k_{0}$ where:

\begin{itemize}
    \item $\mathcal{C}$ is the set of $5$-tuples of non-monochromatic hats. 
    \item $k_{0} = \mathscr{S}_{5,5}$ where: $$\hspace{-20pt} (\mathscr{S}_{5,5}(a_{j}))_{1 \leq j \leq 2^{5}} = (2,3,2,3,5,5,5,5,4,3,2,3,5,5,5,5,1,3,1,3,1,5,1,1,1,3,1,3,1,4,1,5)$$
\end{itemize}

By analogy with $\mathscr{S}_{3}$, we call this strategy $\mathscr{S}_{5}$. Note that when two players use $\mathscr{S}_{5,5}$, they have a probability of winning on five hats equal to $358/2^{10}= 0.349609375$. This probability of winning is computed using a computer by a brute-force algorithm counting the number of winning configurations which equals $358$. It is possible to show \footnote{We present the full computation in Subsection~\ref{subsec:v2(p)}, following the method used in section 5 of \cite{article 3}, adapted to the case
$n=5$.} that the symmetrical strategy $\mathscr{S}_{5}$ has a probability of winning given by:

$$\frac{358-1}{(2^{5}-1)^{2} + 2(2^{5}-1) - 3} = \frac{7}{20}$$

\end{proof}

Notice that it is quite simple to find other strategies that reach 0.35. One method consists of applying the hill-climbing algorithm for $h=5$ in order to obtain a strategy with a probability of victory equal to $0.34960937$. Adapting the previous construction to such a $5$-strategy (i.e., choosing $k_{0}$ to be the $5$-strategy obtained algorithmically) is sufficient. For example, the non-symmetrical $5$-strategy

\begin{center}
$(1, 5, 4, 5, 2, 2, 2, 2, 3, 3, 3, 3, 3, 3, 3, 3, 1, 1, 1, 1, 2, 2, 2, 2, 1, 1, 1, 1, 1, 1, 4, 1),$ \\
$(1, 5, 4, 4, 2, 2, 2, 2, 3, 3, 3, 3, 3, 3, 3, 3, 1, 1, 1, 1, 2, 2, 2, 2, 1, 1, 1, 1, 2, 5, 3, 1)$ 
\end{center}

yields an $\infty$-strategy with probability of victory equal to $0.35$. We emphasize the fact that the hillclimbing algorithm can reach a probability of $0.34960937$ in approximately 10 iterations for $h=5$.  \\

We present a visual representation of $\mathscr{S}_{5}$ along with its associated $5$-strategy. We do not provide an intuitive explanation of the strategy in terms of outcomes due to its complexity. However, it is important to keep in mind that this strategy is built similarly to $\mathscr{S}_{3}$: we start from an efficient $5$-strategy $\mathscr{S}_{5,5}$ and generalize it by 'skipping' monochromatic tuples of size $5$. In the following applications, we will see that, although this strategy achieves the same winning probability as $\mathscr{S}_{3}$, it plays a key role in improving current results.

\begin{figure}[h!]
\centering
\begin{tabular}{cc}

\includegraphics[width=0.3\linewidth]{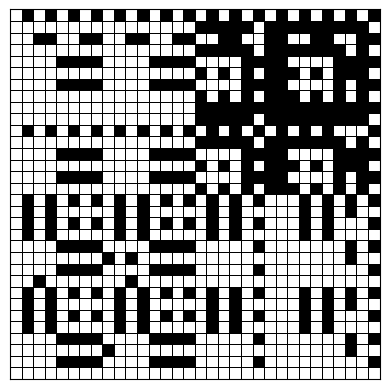} & 
\includegraphics[width=0.3\linewidth]{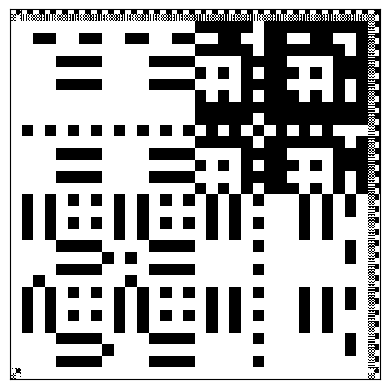} \\

$5$-strategy $\mathscr{S}_{5,5}$ associated to $\mathscr{S}_{5}$  & 
Strategy $\mathscr{S}_{5}$ \\

\end{tabular}
\caption{Left: joint outcome matrix of the $5$-strategy $\mathscr{S}_{5,5}$ for $h=5$ 
(winning probability $358/2^{10} \approx 0.3496$). Right: joint outcome matrix of the 
recursive $\infty$-strategy $\mathscr{S}_5$ (winning probability $7/20 = 0.35$, see 
Theorem~\ref{s5}). As with $\mathscr{S}_3$, the fractal 
structure arises from the recursive skipping of monochromatic $5$-tuples.}
\label{fig:strategies-comparison}
\end{figure}

\subsubsection{Application n°1: A new lower bound for $V_{2,h}$}

In \cite{article 2}, the authors conjecture that the following inequality holds:
\begin{center}
    $V_{2,h} \geq \frac{7}{20} - \frac{C}{r^{h}}$
\end{center}

for some constants $C > 0$ and $r > 1$. Indeed, plotting the values of $V_{2,h}$ suggests a geometric convergence towards $7/20$. However, the ideas presented in \cite{article 2} remained at the level of a conjecture, based on the first few terms of $V_{2,h}$. In this section, we provide a proof of two inequalities of this type respectively based on $\mathscr{S}_{3}$ and $\mathscr{S}_{5}$. We then improve upon it by introducing a larger class of strategies. Of course, such an inequality does not by itself imply convergence to $7/20$. \\

Before turning to the asymptotic behavior of the sequence $(V_{2,h})$, we prove an elementary lemma.

\begin{lemma}
\label{V_2,h != 7/20}
For all $h \geq 1$, we have
\[
\left| \frac{7}{20} - V_{2,h} \right| \geq \frac{1}{5 \cdot 4^h}.
\]
In particular, $V_{2,h} \neq \frac{7}{20}$.
\end{lemma}

\begin{proof}
Let $h \geq 1$. Note that, for the $2$ player game on $h$ hats, there are $4^{h}$ possible configurations. Therefore, the probability of winning using a given strategy in this game is always of the form $p/4^{h}$ where $p$ is the number of winning configurations. However, since the number of $h$-strategies is finite, the supremum defining $V_{2,h}$ is actually a maximum. Equivalently, there exists an $h$-strategy with probability of victory $V_{2,h}$. Therefore, there exists an integer $p \in \mathbb{N}$ such that
\[
V_{2,h} = \frac{p}{4^h}.
\]

Define $\delta_h := \frac{7}{20} - V_{2,h}$. Observe that $5 \cdot 4^h \delta_h \in \mathbb{Z}$ and
\[
5 \cdot 4^h \delta_h = 7 \cdot 4^{h-1} - 5p \not\equiv 0 \pmod{5}.
\]

Hence, in particular, $5 \cdot 4^h \delta_h \neq 0$, which completes the proof.
\end{proof}

This lemma shows—by a purely arithmetic argument—that $h$-strategies never attain the value $7/20$, although they could a priori exceed it. \\

In what follows, we will establish new lower bounds for the sequence $(V_{2,h})_{h \geq 1}$. We begin with two intermediate inequalities.

\begin{proposition}
\label{minoration type K3}
\textbf{(Lower bound of type $\mathscr{S}_{3}$)}\\
For all $h \geq 4$, we have
\[
\frac{7}{20} - V_{2,h} \leq \frac{\frac{7}{20} - V_{2,r}}{16^q},
\]
where $q$ and $r$ are given by the Euclidean division $h = 3q + r$, with $1 \leq r \leq 3$.
\end{proposition}

\begin{proof}
Fix $h \geq 4$, and let $S$ be an optimal strategy for $h$ hats. It is possible to construct a good strategy for $h+3$ hats, denoted $S'$, by proceeding as follows:

\begin{itemize}
    \item If the first three hats are all white or all black, skip them and apply strategy $S$ to the remaining $h$ hats.
    \item Otherwise, on the first three hats, apply the same 3-strategy associated with $\mathscr{S}_{3}$, denoted $\mathscr{S}_{3,3}$.
\end{itemize}

Below is a visualization of $S'$, where $M_S$ denotes the pattern of strategy $S$:

\begin{center}
\includegraphics[width=0.3\textwidth]{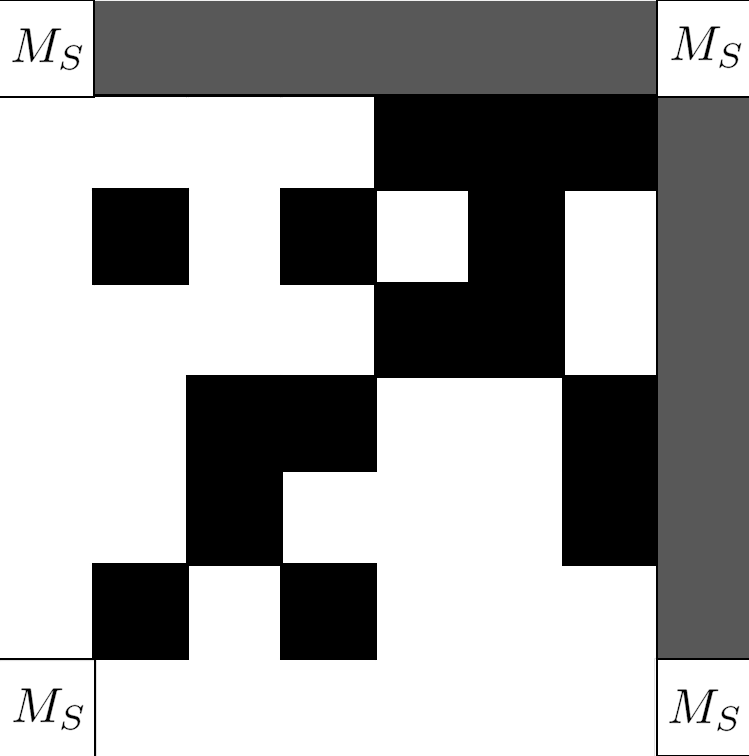}
\captionof{figure}{Visualization of the strategy $S'$ (the shaded area has relative measure $1/2$)}
\end{center}

Recall that the probability that a triplet of hats is monochromatic is $1/4$. There are four possible cases for the two players:

\begin{itemize}
    \item The first three hats of both A and B are monochromatic, with probability $1/16$. Then both players follow strategy $S$ on the remaining $h$ hats; the conditional probability of winning is by construction equal to $V_{2,h}$.
    
    \item The first three hats of both A and B are not monochromatic, with probability $9/16$. Both players apply strategy $\mathscr{S}_{3,3}$ to these hats. The conditional probability of winning is $15/36$.
    
    \item The first three hats of B are monochromatic, but not those of A, with probability $3/16$. The conditional probability of winning is $1/4$. Indeed, there are two equiprobable subcases. If B's three hats are black, B will choose a black hat, and A will choose a black hat with probability $1/2$. In the second subcase, the game is lost.
    
    \item The first three hats of A are monochromatic, but not those of B. This case is symmetric to the previous one and leads to the same result.
\end{itemize}

Therefore, the probability of winning using $S'$ is:
\[
\frac{1}{4} \cdot \frac{1}{4} \cdot V_{2,h} + \frac{3}{4} \cdot \frac{3}{4} \cdot \frac{15}{36} + \frac{1}{4} \cdot \frac{3}{4} \cdot \frac{1}{4} + \frac{1}{4} \cdot \frac{3}{4} \cdot \frac{1}{4}
= \frac{1}{16} V_{2,h} + \frac{21}{64}.
\]

Now, $S'$ is a $(h+3)$-strategy, so $V_{2,h+3} \geq \frac{1}{16} V_{2,h} + \frac{21}{64}$. This yields
\[
\frac{7}{20} - V_{2,h+3} \leq \frac{1}{16} \left( \frac{7}{20} - V_{2,h} \right).
\]

Now write $h = 3q + r$ with $r \in \{1,2,3\}$ and $q \geq 1$. There are two cases to consider:

\begin{itemize}
    \item If $\frac{7}{20} - V_{2,h} \geq 0$, then by monotonicity of $(V_{2,i})$, we have $\frac{7}{20} - V_{2,3l + r} \geq 0$ for all $l \in \{0,1,\dots,q\}$. We may then apply the above inequality recursively and obtain by induction:
\[
\frac{7}{20} - V_{2,h} \leq \frac{\frac{7}{20} - V_{2,r}}{16^q}.
\]

\item 
If $\frac{7}{20} - V_{2,h} < 0$, then $V_{2,h} > 7/20 \geq V_{2,r}$, and in 
particular $\frac{\frac{7}{20} - V_{2,r}}{16^q} \geq 0 \geq \frac{7}{20} - V_{2,h}$, 
so the desired inequality holds trivially. We note however that this case would 
contradict Conjecture~\ref{conj1}, and all numerical evidence suggests it does not occur.
\end{itemize}
\end{proof}

\begin{proposition}
\textbf{(Lower bound of type $\mathscr{S}_{5}$)}\\
For all $h \geq 6$, we have
\[
\frac{7}{20} - V_{2,h} \leq \frac{\frac{7}{20} - V_{2,r}}{256^q},
\]
where $q$ and $r$ are given by the Euclidean division $h = 5q + r$, with $1 \leq r \leq 5$.
\end{proposition}

\begin{proof}
In a similar manner, one can construct from any strategy with $h$ hats a strategy with $h+5$ hats using $\mathscr{S}_{5}$. We denote this new strategy by $S'$. \\

The probability that a fixed player sees the first 5 hats of their partner as monochromatic is $1/16$. As before, we consider the four possible cases. Given that the conditional probability of winning when applying the strategy $\mathscr{S}_{5,5}$—in the case where neither player has their first 5 hats monochromatic—is $109/300$, we obtain the following total probability of winning using $S'$:
\[
\frac{1}{16} \cdot \frac{1}{16} \cdot V_{2,h} + \frac{15}{16} \cdot \frac{15}{16} \cdot \frac{109}{300}
+ \frac{1}{16} \cdot \frac{15}{16} \cdot \frac{1}{4}
+ \frac{1}{16} \cdot \frac{15}{16} \cdot \frac{1}{4}
= \frac{1}{256} V_{2,h} + \frac{357}{1024}.
\]

The remainder of the proof is identical to that of the previous proposition. We thus obtain that for all $h \geq 6$, 
\[
\frac{7}{20} - V_{2,h} \leq \frac{\frac{7}{20} - V_{2,r}}{256^q}.
\]
\end{proof}

We immediately deduce from the lower bound of type $\mathscr{S}_{5}$ the following best asymptotic estimate:

\begin{theorem}
\label{asymptotic_upper_bound_V2h}
In particular, there exists a constant $C > 0$ such that for all $h \geq 1$, 
\[
V_{2,h} \geq \frac{7}{20} - \frac{C}{\left(\sqrt[5]{256}\right)^{h}}.
\]

Note that $\sqrt[5]{256} \approx 3.031\ldots$.\\
\end{theorem}

This result shows in particular that a geometric lower bound on $7/20 - V_{2,h}$ can be established. Furthermore, this inequality allows us to improve the previously presented table:

\begin{longtable}{|c|c|}
    \hline
    \textbf{h} & \textbf{$V_{2,h} \geq$} \\
    \hline
    6 & 0.349609375 \\
    \hline
    7 & 0.349853515625  \\
    \hline
    8 & 0.3499755859375  \\
    \hline
    9 & 0.3499908447265625 \\
    \hline
    10 & 0.34999847412109375  \\
    \hline
    11 & 0.34999847412109375 \\
    \hline
    12 & 0.34999942779541016  \\
    \hline
    13 & 0.34999990463256836 \\
    \hline
    \end{longtable}

This inequality sometimes improves upon the bounds obtained algorithmically via the hill-climbing algorithm. Moreover, if Levine's conjecture holds, then the above upper bound becomes a convergence estimate, yielding the following theorem: \\

\begin{theorem}
\label{thm of convergence}
Assuming that $V_2 = 7/20$, we have the asymptotic estimate
\[
\frac{7}{20} - V_{2,h} = O\left(\frac{1}{256^{h/5}}\right).
\]

Under the same assumption, for all $h \geq 1$, we have
\[
V_{2,h+3} > V_{2,h}.
\]
\end{theorem}

\begin{proof}
Fix $h \geq 1$. The first result follows immediately from the previous theorem and the fact that $V_2 = 7/20$ implies $ V_{2,h} \leq 7/20$. \\

Moreover, we have shown that
\[
V_{2,h+3} \geq \frac{1}{16} V_{2,h} + \frac{21}{64}.
\]
Note that
\[
\frac{1}{16} V_{2,h} + \frac{21}{64} > V_{2,h} \iff \frac{7}{20} > V_{2,h}.
\]

Hence, if $V_2 = 7/20$, then by Lemma \ref{V_2,h != 7/20}, the above inequality holds. This completes the proof.
\end{proof}

The strategy $\mathscr{S}_{5}$ allowed us to obtain a better lower bound for the sequence $(V_{2,h})_{h \geq 1}$ than $\mathscr{S}_{3}$. We now explain how to improve this estimate. First of all, we define a new class of strategies which is a subclass of the recursive strategies. 

\begin{definition}
    Fix $t \geq 3$. We say that an $\infty$-strategy $\mathscr{S}$ is a strategy of type $R_{t}$ if and only if $\mathscr{S}$ is recursive with \begin{itemize}
        \item $\mathcal{C}$ the set of non-monochromatic groups of $t$ consecutive hats
        \item $k_{0}$ is a $t$-strategy which we will often refer to as $\mathscr{S}_{t,t}$
    \end{itemize}
    and such that the probability of victory of $\mathscr{S}$ is exactly $7/20$.
\end{definition}

Note, for example, that $\mathscr{S}_{5}$ and $\mathscr{S}_{3}$ are respectively of type $R_{5}$ and $R_{3}$. It is straightforward to show that, using a strategy of type $R_t$, we obtain a new lower bound on $V_{2,h}$ of the form: 
\[
V_{2,h} \geq \frac{7}{20} - \frac{C_{t}}{\left(4^{1 - \frac{1}{t}}\right)^{h}}.
\]

For some constant $C_{} > 0$ depending on $t$. The key argument is as follows: for such a strategy, the recursive nature implies that the probability of winning in the case where at least one of the two piles is not monochromatic is $X = 7/20$. Indeed, $X$ satisfies the equation
\[
\frac{4}{4^t} \cdot \frac{7}{20} + \left(1 - \frac{4}{4^t}\right) X = \frac{7}{20}.
\]
We do not develop the proof further here, as it is an immediate generalization of the proof of Proposition \ref{minoration type K3}. Let us now study the class of $R_t$-type strategies.\\

Suppose we are given an $R_{t}$-type strategy $\mathscr{S}_{t}$. We naturally denote by $\mathscr{S}_{t,t}$ the associated $t$-strategy. It can then be observed that the conditional probability $X$ mentioned above can be expressed in terms of the winning probability of the strategy $\mathscr{S}_{t,t}$. After some calculations, we obtain the following result.

\begin{lemma}
\label{lemma victory probability K_t,t}
The winning probability of such a strategy $\mathscr{S}_{t,t}$, for $t \geq 3$, is 
\[
\frac{1}{4^t} \left[1 + \frac{7}{5} \left(4^{t-1} - 1\right) \right].
\]
\end{lemma}

\begin{proof}
    Let $p$ denote the winning probability of such a $t$-strategy $\mathscr{S}_{t,t}$, and $\overline{p}$ its conditional winning probability given that both players' stacks are not monochromatic. \\

There are four possible cases:
\begin{itemize}
    \item Both stacks are entirely black, with probability $1/4^t$. In this case, victory is automatic.
    \item One stack is entirely black and the other is non-monochromatic, with probability $2(2^{t}-2)/4^t$. In this case, the winning probability is $1/2$. Indeed, suppose for example that A and B receive these stacks in order. Then A is certain to choose a black hat, and from B we remove two configurations (one winning, one losing) with equal probability, so B will always have a winning probability of $1/2$.
    \item At least one of the two stacks is entirely white. In this case, defeat is automatic.
    \item Both stacks are non-monochromatic, with probability $(2^t-2)^2/4^t$. The winning probability is then $\overline{p}$.    
\end{itemize}

We therefore obtain:

$$p = \frac{2^t-1}{4^t} + \frac{\left(2^t-2\right)^2}{4^{t}}\overline{p}$$

By construction, the associated strategy $\mathscr{S}_t$ coincides with $\mathscr{S}_{t,t}$ on the first $t$ hats of A and B when these are not monochromatic (more precisely, it coincides with the extension of $\mathscr{S}_{t,t}$ to an infinite number of hats, which doesn't change the considered probabilities on the first $t$ hats). By performing the same case analysis for strategy $\mathscr{S}_t$ (using the recursive structure and the fact that the winning probability obtained with $\mathscr{S}_{t}$ is $7/20$), we find:

$$\frac{7}{20} = \frac{1}{4^{t-1}} \frac{7}{20} + \frac{2^t-2}{4^t} + \frac{\left(2^t-2\right)^2}{4^{t}}\overline{p}$$

We obtain the stated result by solving for $p,\overline{p}$.
\end{proof}

However, $t$-strategies have a winning probability that is an integer multiple of $1/4^t$. Moreover, $4^{t-1} - 1 \equiv 0 \pmod{5}$ if and only if $t$ is odd. Thus, we obtain the following theorem:

\begin{theorem}
There exists no $R_{t}$-type strategy for even integers $t \geq 3$.
\end{theorem}

At this point, it is natural to ask whether $R_t$-type strategies exist for infinitely many odd $t$'s. Indeed, as we saw previously, having an infinite amount of such strategies would yield a geometric convergence rate of ratio $1/4^{1-\varepsilon}$ for any $1 > \varepsilon > 0$. In fact, we will now show that $R_t$-type strategies exist for every odd $t \geq 3$. We are grateful to Henry Swanson \cite{article 7} for sharing key insights regarding the construction of $R_t$-type strategies following an earlier preprint version of this work. We start with a few important definitions:

\begin{definition}[Standardized Strategy]
A finite strategy $\mathscr{S}$ is said to be standardized if, when a player observes a monochromatic block of hats on their partner's head, the strategy returns the same fixed index. In other words, the player's choice of index is invariant to the specific color of the monochromatic block they are looking at.
\label{def:standardized}
\end{definition}

\noindent We now define the main pillar of the proof. We explain how to construct a $(N+M-1)$-strategy using a $N$-strategy and a $M$-strategy which are standardized.

\begin{definition}[Combining Construction]
Let $\mathscr{S}$ be a standardized $N$-strategy and $\mathscr{T}$ a standardized $M$-strategy. We define the composed $(N+M-1)$-strategy $\mathscr{S} \star \mathscr{T}$ as follows:
\begin{enumerate}
    \item If the first $N$ hats observed are non-monochromatic, the player applies $\mathscr{S}$ to these hats and returns the resulting index $i \in \{1, \dots, N\}$.
    \item If the first $N$ hats are monochromatic, the player applies $\mathscr{T}$ to the hats at positions $N+1, \dots, N+M$ and returns the resulting index offset by $N$, i.e., $i \in \{N+1, \dots, N+M\}$.
\end{enumerate}
\end{definition}

\noindent We now wish to assess the probability of victory of these "combined" strategies.

\begin{lemma}[Combining Formula]
\label{combining formula}
Let $w(\mathscr{S})$ and $w(\mathscr{T})$ denote the winning probabilities of the standardized strategies $\mathscr{S}$ and $\mathscr{T}$ in the $2$-player game. Suppose that $\mathscr{S}$ is an $N$-strategy and $\mathscr{T}$ an $M$-strategy. The composed strategy $\mathscr{S} \star \mathscr{T}$ has a winning probability of:
\begin{equation}
    w(\mathscr{S} \star \mathscr{T}) = w(\mathscr{S}) + \frac{4w(\mathscr{T}) - 1}{4^N}
\end{equation}
\end{lemma}

\begin{proof}
Let Alice and Bob play with strategy $\mathscr{S} \star \mathscr{T}$. We partition the configurations of the first $N$ hats on each player's head into three cases:
\begin{itemize}
    \item \textbf{Both non-monochromatic:} This occurs with probability $\frac{(2^N-2)^2}{4^N}$. Both players apply $\mathscr{S}$. Let $p^*$ be the conditional winning probability in this case.
    \item \textbf{Both monochromatic:} This occurs with probability $\frac{4}{4^N}$ (since each player independently sees a monochromatic block with probability $\frac{2}{2^N}$). In this case, both players apply $\mathscr{T}$ to the subsequent $M$ hats. Since these hats are independent of the first $N$, the conditional winning probability is exactly $w(\mathscr{T})$.
    \item \textbf{Exactly one monochromatic:} This occurs with probability $\frac{4(2^N-2)}{4^N}$. Suppose Alice's first $N$ hats are monochromatic. She picks an index $i_A > N$ based on $\mathscr{T}$. Bob, seeing Alice's non-monochromatic block, picks an index $i_B \leq N$ based on $\mathscr{S}$. Because Bob's $i_B$-th hat belongs to a monochromatic block on his head, it is black with probability $1/2$. Alice's $i_A$-th hat is also black with probability $1/2$ and is independent of $i_B$. Thus, the conditional winning probability is $\frac{1}{2} \times \frac{1}{2} = \frac{1}{4}$.
\end{itemize}

By applying the same case analysis to the isolated strategy $\mathscr{S}$, and noting that for a standardized strategy, the win probability in any monochromatic case is $1/4$, we have:
\[ w(\mathscr{S}) = p^* \frac{(2^N-2)^2}{4^N} + \frac{1}{4} \frac{4(2^N-2)}{4^N} + \frac{1}{4} \frac{4}{4^N} = p^* \frac{(2^N-2)^2}{4^N} + \frac{2^N-1}{4^N} \]
Substituting $p^*$ into the expression for $w(\mathscr{S} \star \mathscr{T})$:
\[ w(\mathscr{S} \star \mathscr{T}) = \left( w(\mathscr{S}) - \frac{2^N-1}{4^N} \right) + \frac{1}{4} \frac{4(2^N-2)}{4^N} + w(\mathscr{T}) \frac{4}{4^N} \]
\[ w(\mathscr{S} \star \mathscr{T}) = w(\mathscr{S}) + \frac{-(2^N-1) + (2^N-2) + 4w(\mathscr{T})}{4^N} = w(\mathscr{S}) + \frac{4w(\mathscr{T}) - 1}{4^N} \]
\end{proof}

\noindent We may now use these tools to obtain the final required result:

\begin{theorem}
\label{thm:Rt_exists}
For all odd $t \geqslant 3$, there exists an $R_t$-type strategy.
\end{theorem}

\begin{proof}
We proceed by strong induction on odd $t$ to show two results simultaneously:
\begin{enumerate}
    \renewcommand{\theenumi}{\roman{enumi}}
    \renewcommand{\labelenumi}{(\theenumi)}
    \item there exists a standardized strategy $\mathscr{S}_{t,t}$ achieving a winning 
    probability of
    \begin{equation}
        \label{eq:wt}
        w(\mathscr{S}_{t,t}) = \frac{7}{20} - \frac{2}{5 \cdot 4^t}
    \end{equation}
   
    which corresponds to the probability required by Lemma~\ref{lemma 
    victory probability K_t,t};
    \item the recursive infinite strategy built upon $\mathscr{S}_{t,t}$ is an 
    $R_t$-type strategy.
\end{enumerate}The simultaneous induction on standardizedness is necessary to apply the 
Combining Construction at each step.

\paragraph{Base case ($t = 3$).}
The strategy $\mathscr{S}_{3,3}$ is defined by
\[
(k_1(a_j))_{1 \leqslant j \leqslant 8} = (1,3,2,2,1,3,1,1), \quad
(k_2(a_j))_{1 \leqslant j \leqslant 8} = (1,3,2,3,1,1,2,1).
\]
The all-white stack $a_1 = (0,0,0)$ and the all-black stack $a_8 = (1,1,1)$ 
both return index $1$. Hence $\mathscr{S}_{3,3}$ is standardized. Its winning 
probability is $w(\mathscr{S}_{3,3}) = \frac{22}{64} = \frac{11}{32}$, and 
one checks that
\[
\frac{7}{20} - \frac{2}{5 \cdot 4^3} = \frac{11}{32},
\]
so \eqref{eq:wt} holds. The corresponding recursive strategy $\mathscr{S}_3$ 
is an $R_3$-type strategy, establishing the base case.

\paragraph{Inductive step.}
Let $t \geqslant 5$ be an odd integer, and assume the inductive hypothesis 
holds for $t-2$: there exists a standardized strategy $\mathscr{S}_{t-2,t-2}$ 
satisfying $w(\mathscr{S}_{t-2,t-2}) = \frac{7}{20} - \frac{2}{5 \cdot 4^{t-2}}$.
We construct $\mathscr{S}_{t,t} := \mathscr{S}_{t-2,t-2} \star \mathscr{S}_{3,3}$.

\medskip
\noindent\textit{Standardizedness of $\mathscr{S}_{t,t}$.}
If a player observes a monochromatic $t$-tuple, then in particular the first 
$t-2$ hats are monochromatic. By the Combining Construction, the player then 
applies $\mathscr{S}_{3,3}$ to the remaining hats. Since the full $t$-tuple 
is monochromatic, those remaining hats are also monochromatic, and since 
$\mathscr{S}_{3,3}$ is standardized, it returns the fixed index $1$ regardless 
of color. Hence $\mathscr{S}_{t,t}$ returns index $(t-2)+1 = t-1$ for any 
monochromatic $t$-tuple, independently of color. Thus $\mathscr{S}_{t,t}$ is 
standardized.

\medskip
\noindent\textit{Winning probability.}
Since both $\mathscr{S}_{t-2,t-2}$ and $\mathscr{S}_{3,3}$ are standardized, 
the Combining Formula (Lemma~\ref{combining formula}) applies with $N = t-2$:
\[
w(\mathscr{S}_{t,t}) = w(\mathscr{S}_{t-2,t-2}) + 
\frac{4\,w(\mathscr{S}_{3,3}) - 1}{4^{t-2}}.
\]
Substituting $w(\mathscr{S}_{3,3}) = \frac{11}{32}$ gives 
$4 \cdot \frac{11}{32} - 1 = \frac{3}{8}$, and using the inductive hypothesis:
\[
w(\mathscr{S}_{t,t}) = \left(\frac{7}{20} - \frac{2}{5 \cdot 4^{t-2}}\right) 
+ \frac{3/8}{4^{t-2}} = \frac{7}{20} - \frac{1}{40 \cdot 4^{t-2}} 
= \frac{7}{20} - \frac{2}{5 \cdot 4^t},
\]
which matches \ref{eq:wt}. Consequently, the recursive infinite strategy 
built upon $\mathscr{S}_{t,t}$ by taking $\mathcal{C}$ to be the set of 
non-monochromatic $t$-tuples achieves the limit winning probability of 
$7/20$, making it an $R_t$-type strategy.
\end{proof}

\noindent As explained earlier, we may now conclude on the speed of convergence of $V_{2,h}$. \\

\begin{theorem}[Rate of convergence of $(V_{2,h})$] \label{final_rate}For all $\varepsilon >0$, there exists a constant $C_{\varepsilon} > 0$ such that: \begin{center}
    $\forall h \geq 1, V_{2,h} \geq \frac{7}{20} - \frac{C_{\varepsilon}}{\left(4^{1-\varepsilon} \right)^h}$
\end{center}
    
\end{theorem}

\begin{proof}
Given $\varepsilon > 0$, choose an odd integer $t \geqslant 3$ such that $\frac{1}{t} \leqslant \varepsilon$, 
which is possible by Theorem~\ref{thm:Rt_exists}. By the argument above, 
the $R_t$-type strategy yields, for all $h \geqslant 1$,
\[
V_{2,h} \geqslant \frac{7}{20} - \frac{C_t}{\left(4^{1-1/t}\right)^h},
\]
for some constant $C_t > 0$. Since $4^{1-1/t} \geqslant 4^{1-\varepsilon}$, 
setting $C_\varepsilon := C_t$ yields the result.
\end{proof}

\subsubsection{Application n°2: Improvement of a Known Bound, variant $p \neq 1/2$}
\label{subsec:v2(p)}
Previous studies ask what happens when the probability distribution in Levine's hat problem is modified. Specifically, one may assign to each hat a probability \( p \in ]0,1[ \) of being black (whereas in the original problem, \( p = 1/2 \)). This generalization of the problem is one of the approaches introduced in reference \cite{article 3}. In particular, authors show that efficient strategies are not necessarily the same for different values of \( p \). \\

Our formalization introduced in previous sections can be extended to the case \( p \in ]0,1[ \). This extension corresponds to a distortion of the original measure according to \( \mu(C(\varepsilon)) = p^i (1-p)^{h-i} \), where \( i \) denotes the number of 1's appearing in \( \varepsilon \in \Omega_h \). In particular, a graphical representation associated with this deformed measure can be provided. We consider the example of the first black hat strategy, depicted below for various values of \( p \). In this example, we observe geometrically that the winning probability (still corresponding to the area of the associated domain) increases with \( p \), due to the change in the measure.

\begin{figure}[h!]
    \centering
    \begin{tabular}{ccc}
        \includegraphics[width=0.3 \linewidth]{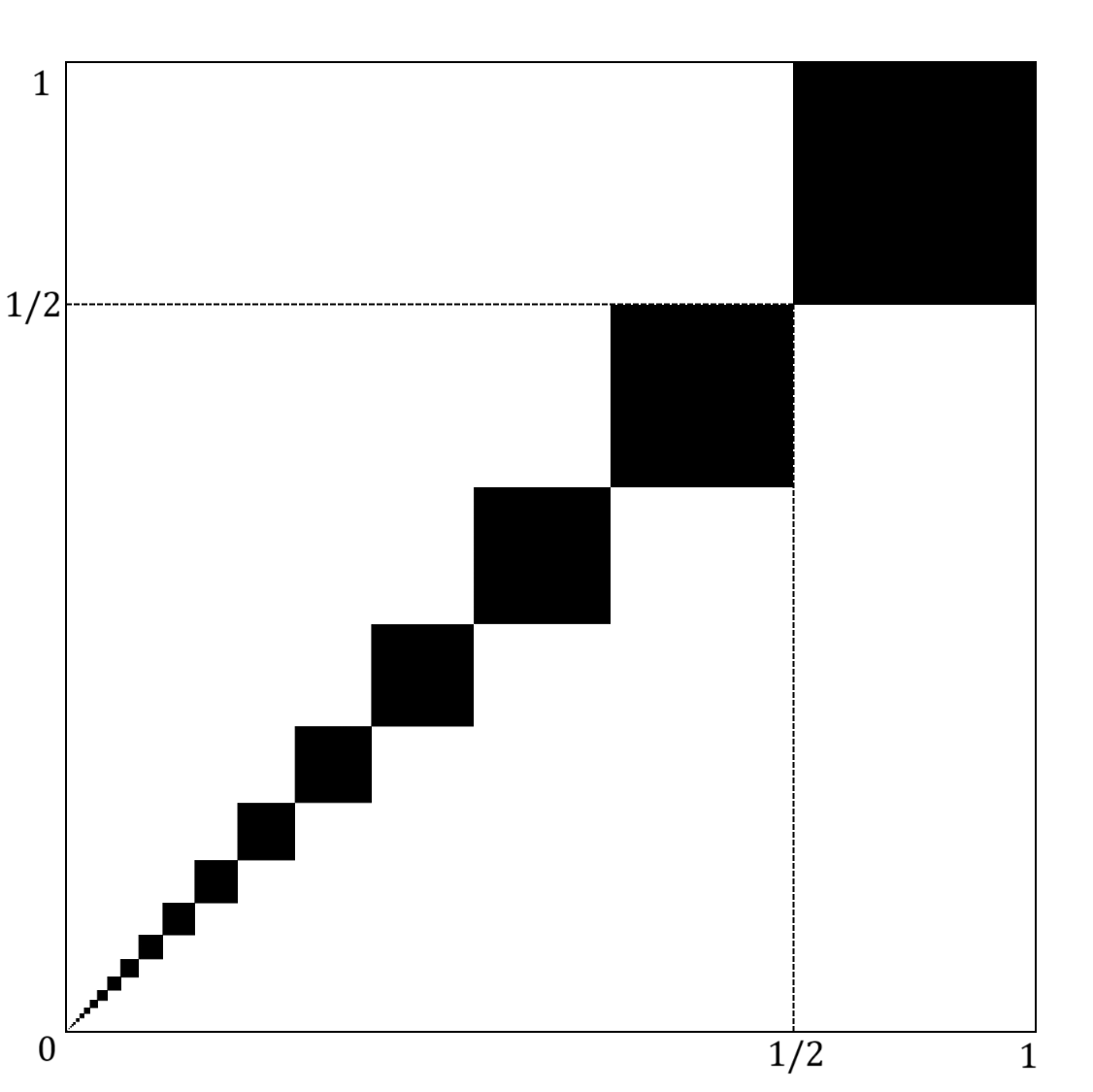} & 
        \includegraphics[width=0.3 \linewidth]{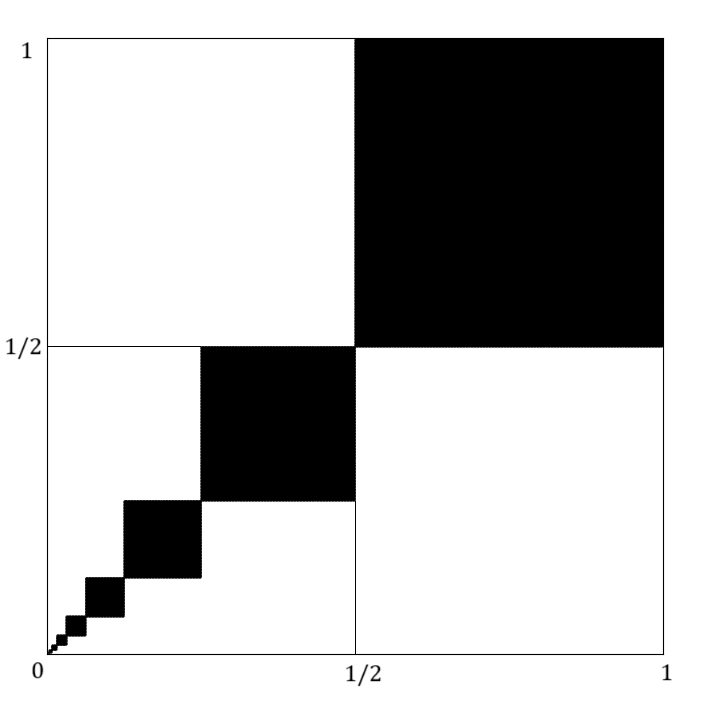} &
        \includegraphics[width=0.3 \linewidth]{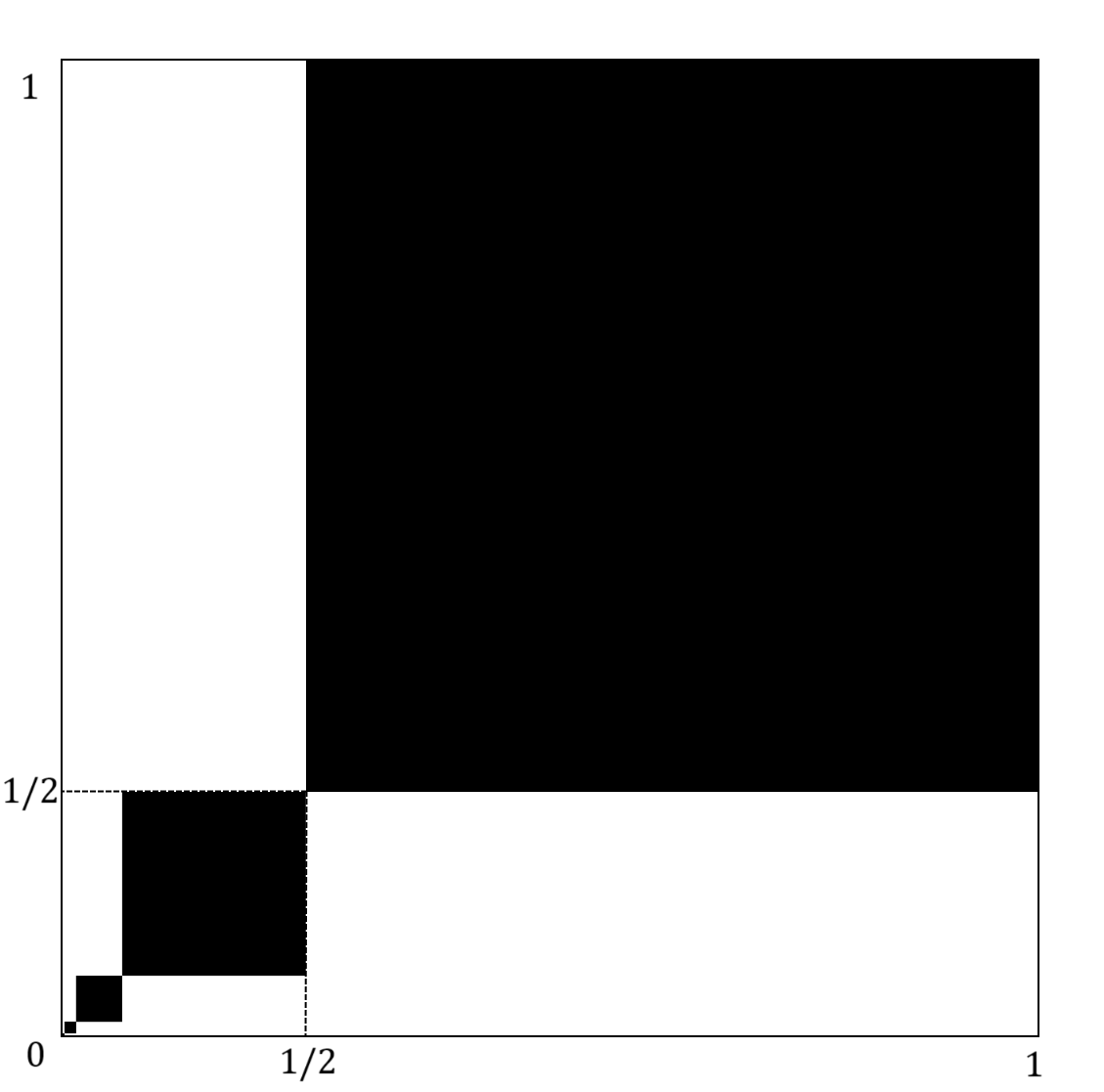} \\
        $p=1/4$ & $p=1/2$ &
$p=3/4$ \\
    \end{tabular}
    \caption{Visualization of the first black hat strategy ($n=2$) for various values of $p$}

    %\label{fig:fbh_3h}
\end{figure}

\noindent The article \cite{article 3} investigates the optimal winning probability in this game, denoted by $V_2(p)$ (we also denote by $V_2(p,\mathscr{S})$ the winning probability when using a fixed strategy $\mathscr{S}$ in this setting). In particular, it establishes the following lower bounds, by reusing good strategies known from the case $p=1/2$.

\begin{proposition} (Joe Buhler \textit{et al.}, theorem 3 in \cite{article 3})

The following holds

\begin{itemize}
    \item For $p \leq 1/2$, $V_2(p) \geq \underbrace{\frac{p + p^2 + p^{3} + 3p^{4} - 3p^{5} + p^{6}}{2+p+p^2+p^3-p^4}}_{:=U_1(p)}$
    \item For $p \geq 1/2$, $V_2(p) \geq \underbrace{\frac{p + 5p^2 - 10p^{3} + 10p^{4} - 5p^{5} + p^{6}}{4-2p-2p^2+3p^3-p^4}}_{:=U_2(p)}$
\end{itemize}

\end{proposition}

We now aim to exhibit a sharper lower bound for $V_{2}(p)$ than those previously known. Moreover, we provide a general method to further improve such bounds. Once again, the starting point is the strategy $\mathscr{S}_{5}$ that we introduced. This strategy allows us to establish a new lower bound $U_{3}(p)$, which is strictly better on a non-empty domain.

\begin{theorem}
\label{u3(p)}
For all $0 < p < 0.312\ldots$, we have $V_{2}(p) \geq U_{3}(p) > \max(U_{1}(p),U_{2}(p))$ \\

where 
\[
U_{3}(p) = \frac{5p-20p^{2} +51p^{3}-82p^{4} + 85p^{5} - 52p^{6} + 10p^{7} +10p^{8} - 7p^{9}}{10-45p + 120p^{2} - 210p^{3} + 250p^{4} - 200p^{5} + 100p^{6} - 25p^{7}}.
\]
\end{theorem}

\begin{proof}
Let $p \in ]0,1[$. Set $q=1-p$. Consider the 5-strategy associated with $\mathscr{S}_{5}$, still denoted $\mathscr{S}_{5,5}$. We study the game with an infinite number of hats. There are $4$ distinct cases:

\begin{itemize}
    \item The first 5 hats of A and B are not monochromatic. The probability that this occurs and the players win (by applying $\mathscr{S}_{5,5}$) is
    $$3p^{10} - 10p^{9} + 30p^{8} - 62p^{7} + 85p^{6} - 82p^{5} + 51p^{4} - 20p^{3} + 5p^{2}.$$ We derived this exact quantity using a computer by enumerating every possible such stack combinations and summing their probability contributions.  
    \item The first 5 hats of A are black, and those of B are not monochromatic with probability $p^{5}(1-p^{5}-q^{5})$. The probability of winning in this case is $p^{6}(1-p^{5}-q^{5})$ because A will automatically choose a black hat and B will choose one with probability $p$. The symmetric situation in A/B yields the same winning probability.    
    \item The first 5 hats of A are white, and those of B are not monochromatic. Then the players lose automatically. The symmetric situation in A/B yields the same winning probability.    
    \item The first 5 hats of A and B are monochromatic with probability $(p^{5}+q^{5})^{2}$. By definition of the strategy $\mathscr{S}_{5}$, they then win with probability $(p^{5}+q^{5})^{2}V_2(p,\mathscr{S}_{5})$.
\end{itemize} 

Thus, we have 

\begin{align*}
\begin{split}
    V_{2}(p,\mathscr{S}_{5}) = {} & 3p^{10} - 10p^{9} + 30p^{8} - 62p^{7} + 85p^{6} - 82p^{5}  + 51p^{4} - 20p^{3} \\ 
    &+ 5p^{2} + 2p^{6}(1-p^{5}-q^{5}) + (p^{5}+q^{5})^{2}V_2(p,\mathscr{S}_{5})
\end{split}
\end{align*}

Hence the result upon solving for $V_{2}(p, \mathscr{S}_{5})$ and expanding. The relevant range for $p$ has been verified by a computer.
    
\end{proof}

\begin{figure}[h!]
     \centering
     \includegraphics[width=1\linewidth]{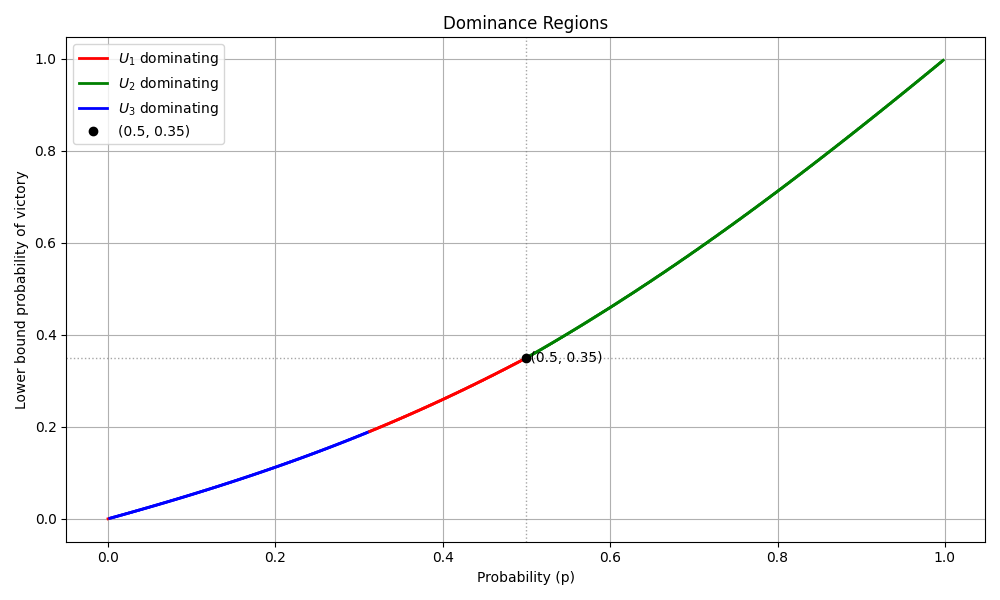}
     \caption{Plot of $\max(U_{1}(p),U_{2}(p),U_{3}(p))$}
     \label{fig:enter-label}
 \end{figure}

It is likely that this bound can be improved even further by finding other recursive strategies based on $h$-strategies, just as we built $\mathscr{S}_{5}$. In fact, $U_{1}$ and $U_{2}$ were obtained using distinct $R_3$-type strategies. Hence it is likely that we can obtain a better bound on a different interval using only strategies from $R_5$. We leave this direction to the interested reader. \\

Furthermore, our use of a recursive strategy of order $5$ provides a negative answer to an expectation raised in \cite{article 3}. Indeed, the authors of \cite{article 3} believed that recursive strategies of order $t > 3$ would not be better than those of order $3$. While $\mathscr{S}_5$ only improve $V_2(p)$ for a certain range relatively far from $p=0.5$, it has also improved the convergence rate of $\left( V_{2,h} \right)_{h \geq 1}$. In fact, these results seem to show that strategies improve as $t$ grows. \\

The constructions of Section~\ref{section algo} suggest that the difficulty of the problem is intimately tied to the discrete binary structure of hat stacks. The continuous generalization introduced here is designed to isolate this obstruction precisely.

\section{A New Generalization: The continuous Levine Game}

\subsection{Statement of the continuous Levine Game}

The study of Levine's hat problem has naturally given rise to a variety of 
generalizations, each designed to isolate the structural features that make 
the original problem difficult. Buhler et al.\ \cite{article 3} study two 
such variants: the biased coin setting where each hat is black with 
probability $p \neq 1/2$, which reveals how strategies conjectured to be optimal shift 
with $p$; and the restricted variant where each player must 
point to their first black hat, which admits a clean analysis and serves 
as a baseline for understanding the general case. Aubrun et al.\ 
\cite{article 2} consider a different relaxation in which the team wins if 
all players simultaneously point to a hat of the same color, whether black 
or white, a condition which grants the players strictly more winning 
configurations to exploit. \\

We introduce here a new generalization in this spirit: the continuous 
Levine game, in which players are granted an uncountably infinite set of 
hats derived from their original stack. We show that the optimal winning 
probability in this setting equals exactly $1/2$ for all $n \geq 2$. Each 
of the generalizations above illuminates, by contrast, which constraint is 
doing the essential work in the original problem. Our result identifies the 
discrete binary structure of the hat stacks as what keeps $V_n$ strictly 
below $1/2$: by relaxing it, the players gain enough freedom to perfectly 
correlate their strategies, and the problem reduces from a combinatorial 
optimization to a tractable analytic one.

At first, each of the $n \geq 2$ players is assigned a countably infinite sequence of black $(1)$ or white $(0)$ hats, according to independent Bernoulli($1/2$) distributions. Then, an uncountable set of additional hats is assigned to each player, determined from their initial sequence in the following way. We denote by
\[
\left(X^{(2^k)}_i\right)_{k \geq 1}
\]
the initial countable set of hats assigned to player $i$, and we consider the real number $x_i$ encoded in base $2$ by this sequence:

$$x_i := \sum_{k \geq 1} \frac{X_i^{(2^{k})}}{2^{k}}$$

Player $i$ is then assigned a larger set of hats indexed by real numbers $a \in \mathbb{R}_+$, and defined by

$$X_i^{(a)} := \lfloor ax_{i}\rfloor \text{ mod } 2$$

As in the original Levine game, each player $i$ can observe the hat sets of the other players, but not their own. Each player must then \textit{choose}\footnote{ using a Lebesgue-measurable strategy, as we explain further below.} a real number $a_i \geq 0$. The team wins if and only if, for every $i$, 
\[X_i^{(a_i)} = 1.\]
Of course, the players may agree on a strategy before the game begins, but they cannot communicate once the game starts. \\

Note that the construction of the continuous hat sets may appear ambiguous for indices of the form $a = 2^\ell$ with $\ell \geq 1$. In fact, the two definitions coincide almost surely. Indeed, for every $\ell \geq 1$ and for every family
\[
\left( X^{(2^k)} \right)_{k \geq 1} \in \{0,1\}^{\mathbb{N}^*}
\]
which is non-stationary at $1$ from some rank onward, we write

$$ \left \lfloor 2^\ell \sum_{k \geq 1} \frac{X^{(2^{k})}}{2^{k}} \right \rfloor = \left \lfloor \sum_{k \geq 1} X^{(2^{k})}2^{\ell-k}  \right \rfloor = \sum_{k = 1}^\ell X^{(2^{k})}2^{\ell-k} = X^{(2^{\ell})} \text{ mod } 2$$

At first glance, it might not seem clear why this continuous extension of the game is natural. We now explain this point more clearly. In the original game, the players are allowed to look at hats indexed by positive integers. However, we saw that picking the hat of index $k$ was strictly equivalent to choosing the bit of index $k$ of a real number in $[0,1[$. We extend the formula that yields the bit of index $k \in  \N^{*}$ to allow us to extract a bit of "index" $a \geq 0$. This extends the sets of hats by filling the gaps and allowing the players to choose their respective indices within a continuous set (namely $\R^{+}$).  \\

We emphasize the key idea of this generalization: it is a continuous extension of the set of hats, which fundamentally provides greater freedom of choice in the game. The players are free to ignore the hats whose indices are not integer powers of $2$ (which is equivalent to considering Levine's original problem). Thus, the optimal probability of winning in the continuous game is greater than $V_n$. \\

Another property that makes this construction natural is the following. The probability that a hat with fixed index $a \geq 0$ is black is by definition

$$ \int_0^1 \textbf{1} \left( \lfloor ax\rfloor \equiv 1 \text{ mod } 2\right) \D x$$

which can be rewritten with the previous $\varepsilon$ function as

$$ \int_0^1 \varepsilon(ax) \D x.$$

However, the following result holds.

\begin{lemma}
\label{integral epsilon <= 1/2}
For all $a \geq 0$, 
$$\int_{0}^{1} \varepsilon(ax)\D x \leq \frac{1}{2}$$
\end{lemma}

\begin{proof}
The case $a = 0$ is straightforward. We therefore fix $a > 0$ and observe that:
\[
\int_{0}^{1} \varepsilon(ax)\, \mathrm{d}x = \frac{1}{a} \int_{0}^{a} \varepsilon(x)\, \mathrm{d}x.
\]
Define $\gamma(a) := \frac{a}{2} - \int_{0}^{a} \varepsilon(x)\, \mathrm{d}x$. We note that $\gamma$ is $2$-periodic. Indeed:
\begin{align*}
    \gamma(a+2) &= 1 + \frac{a}{2} - \int_{0}^{a} \varepsilon(x)\, \mathrm{d}x - \int_{a}^{a+2} \varepsilon(x)\, \mathrm{d}x \\
                &= 1 + \gamma(a) - \int_{0}^{2} \varepsilon(x)\, \mathrm{d}x \hspace{5pt} \text{(by the $2$-periodicity of } \varepsilon) \\
                &= \gamma(a).
\end{align*}

Moreover, if $0 \leq a \leq 1$, then $\gamma(a) = \frac{a}{2}$. If $1 \leq a \leq 2$, then $\gamma(a) = 1 - \frac{a}{2}$. This means that the following holds:
\[
\forall a \geq 0, \quad \gamma(a) \geq 0.
\]

In particular, we deduce that $\int_{0}^{1} \varepsilon(ax)\, \mathrm{d}x \leq \frac{1}{2}$ for all $a > 0$.
\end{proof}

Thus, no player can select a black hat from their own stack with a probability greater than $1/2$. Equivalently, the additional freedom granted to the players remains comparable to that of the original game. \\

Finally, although this variant is based on additional freedom of choice, it does not provide the players with more information. Indeed, the new hats whose indices are not integer powers of $2$ are entirely determined by the initial sequences of hats. Thus, this game can yield interesting insights into the original problem.

\subsection{Formulation of Imaginary Strategies}

We now focus on the definition of strategies for this game, which we refer to as \textbf{imaginary strategies}. Since the final set of hats is fully determined by the same set of hats as in the original game (and with the same distribution), one should use Lebesgue-measurable functions as strategies. In this way, the class of imaginary strategies generalizes that of the original strategies.  \\

\begin{definition}
We denote by $\widehat{\mathcal{M}}_{n-1}$ the set of measurable, \textbf{non-negative} functions on $[0,1]^{n-1}$ for $n \geq 2$, which we call \emph{imaginary} strategies for $n$ players. \\

We then define:
\[
W_n := \sup_{f_1, \dots, f_n \in \widehat{\mathcal{M}}_{n-1}} \int_{[0,1]^n} \prod_{i=1}^n \varepsilon \bigl(f_i(x_{-i}) x_i \bigr) \, \mathrm{d}\Lambda^n,
\]

where $\varepsilon: x \longmapsto \lfloor x \rfloor - 2\lfloor x/2\rfloor$. It follows that $W_n \geq V_n$.
  
\end{definition}

\begin{theorem}
\label{thm W_n >= V_n}
The optimal winning probability in the $n$-player continuous game is $W_n$.
\end{theorem}

\begin{proof}
    Let $f_1, \dots, f_n \in \widehat{\mathcal{M}}_{n-1}$ be the imaginary strategies of the players. The strategy $f_i$ returns the real index chosen by player $i$ based on the hats of the other players. We denote by $X_i$ the stack of player $i$ and by $X_i^{(a)}$ their hat with index $a \geq 0$. We also denote by $X_{-i}$ the vector formed by the $X_j$ for $j \neq i$. The probability of the players' victory is then:
    
    \begin{align*}
        \Proba \left(\forall 1 \leq i \leq n, X^{(f_i(X_{-i}))}_{i} = 1 \right) &= \int_{[0,1]^{n}} \prod_{i=1}^{n}  \mathbf{1}\left(\lfloor f_{i}(x_{-i})x_{i}) \rfloor \equiv 1 \text{ mod } 2 \right) \mathrm{d} \Lambda^{n} \\ 
        &= \int_{[0,1]^{n}} \prod_{i=1}^{n}\varepsilon(f_{i}(x_{-i})x_{i}) \mathrm{d} \Lambda^{n}
    \end{align*}

Indeed, for all $a,x \geq 0$, we have : $\varepsilon(ax)=1$ if and only if $\lfloor ax \rfloor = 1$ mod 2, hence the result.
\end{proof}

At this point, we may obtain valuable information about $V_{n}$. In fact, this continuous hat game can be solved.

\subsection{Optimality in Levine's continuous hat game}

We already saw that the probability of finding a black hat for a single player is always smaller than $1/2$. Hence the probability of a collective victory with any strategy in the continuous hat game is also smaller than $1/2$. In fact, this bound is sharp. To prove this, we will first need the following elementary lemma:

\begin{lemma}
\label{epsilonto1/2}
It holds that
\[ \lim_{a \to \infty} \int_0^1 \varepsilon(ax) \, dx = \frac{1}{2}. \]
\end{lemma}

\begin{proof}
We first observe that $\varepsilon(x)$ is periodic with period $T = 2$. Indeed, for any $x \in \mathbb{R}$:
\[ \varepsilon(x+2) = \lfloor x+2 \rfloor - 2\left\lfloor \frac{x+2}{2} \right\rfloor = (\lfloor x \rfloor + 2) - 2(\lfloor x/2 \rfloor + 1) = \varepsilon(x). \]
Furthermore, $\varepsilon(x)$ takes values in $\{0, 1\}$. Specifically, $\varepsilon(x) = 0$ for $x \in [0, 1)$ and $\varepsilon(x) = 1$ for $x \in [1, 2)$. The integral over one period is thus $\int_0^2 \varepsilon(t) \, dt = 1$.

Consider the integral $I(a) = \int_0^1 \varepsilon(ax) \, dx$. Applying the substitution $t = ax$, we obtain:
\[ I(a) = \frac{1}{a} \int_0^a \varepsilon(t) \, dt. \]
Let $N = \lfloor a/2 \rfloor$ be the number of full periods in the interval $[0, a]$. We decompose the integral into the sum of $N$ full periods and a remainder $R$:
\[ I(a) = \frac{1}{a} \left( \sum_{k=0}^{N-1} \int_{2k}^{2(k+1)} \varepsilon(t) \, dt + \int_{2N}^a \varepsilon(t) \, dt \right) = \frac{N}{a} + \frac{R}{a}, \]
where $R = \int_{2N}^a \varepsilon(t) \, dt$. Given that $0 \leq \varepsilon(t) \leq 1$ and the remainder interval length satisfies $a - 2N < 2$, it follows immediately that $0 \leq R < 2$.

As $a \to \infty$, we note that $\frac{N}{a} = \frac{\lfloor a/2 \rfloor}{a} \to \frac{1}{2}$. Simultaneously, since $R$ is bounded, the term $\frac{R}{a}$ vanishes. Thus,
\[ \lim_{a \to \infty} I(a) = \frac{1}{2} + 0 = \frac{1}{2}. \]
\end{proof}

We may now prove the main theorem of this section:

\begin{theorem}
    For all $n\geq 2$, $W_{n} = 1/2$.
\end{theorem}

\begin{proof}

We define a natural sequence of imaginary strategies for $n$ players that provide this result. Consider for each $m \geq 1$,
$$f_{1,m} = \dots = f_{n,m} = m\pi_{n-1}$$
where $\pi_{n-1}: (x_{1}, \dots, x_{n-1}) \mapsto x_{1}\dots x_{n-1}$. \\

Notice that all of the terms of the form $\varepsilon(f_{i,m}(x_{-i})x_{i}), 1 \leq i \leq n$ are equal to $\varepsilon(mx_{1}\dots x_{n}) \in \{0,1\}$. This means that the factors in the product appearing in $W_{n}$'s integral expression are perfectly correlated (i.e the players always win collectively and always lose collectively) and yield a probability of victory equal to: \begin{center}
    $p_m = \int_{[0,1]^{n}} \varepsilon(mx_{1}\dots x_{n}) \text{d}\Lambda^{n}$
\end{center}

We present the visualization of these strategies for $n=2$ below. It is possible to show that $p_{m}$ tends to $1/2$ as $m$ goes to infinity. \\

\begin{figure}[h!]
    \centering
    \begin{tabular}{ccc}
        \includegraphics[width=0.3 \linewidth]{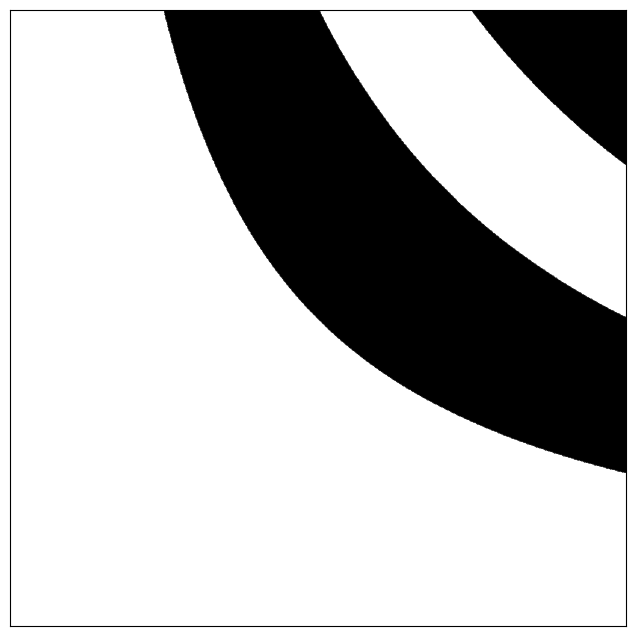} & 
        \includegraphics[width=0.3 \linewidth]{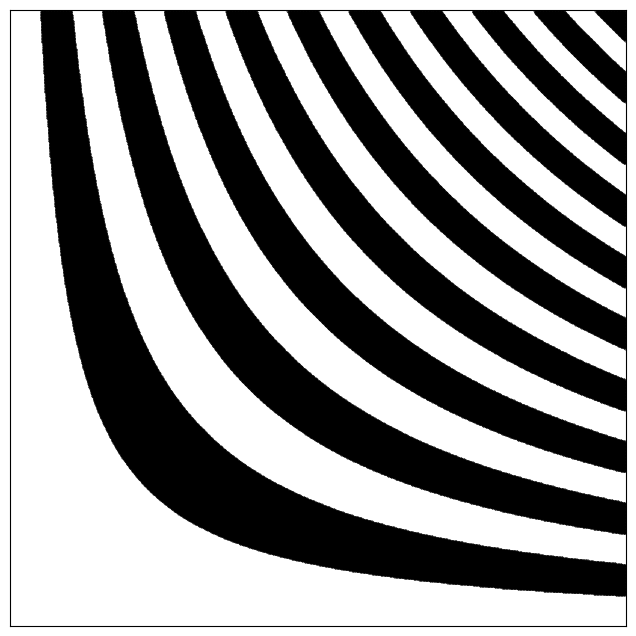} &
        \includegraphics[width=0.3 \linewidth]{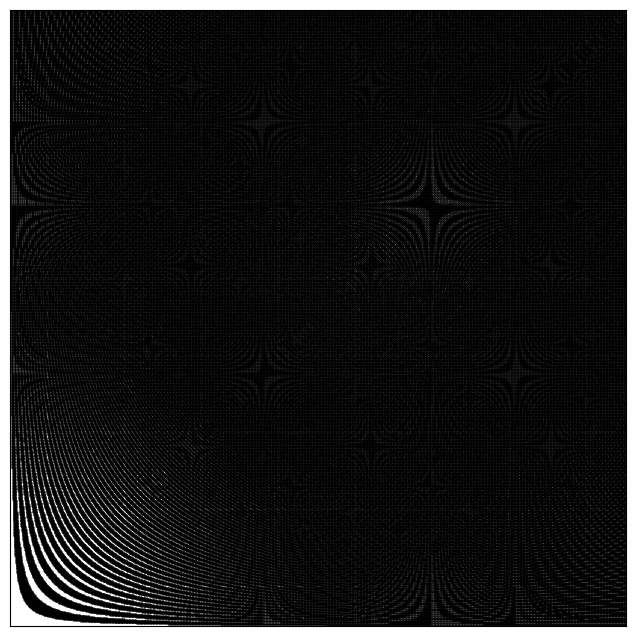} \\
        $m=4$, $p_m \simeq 0.28$ & $m=20$, $p_m \simeq 0.44$ &
        $m=10^3$, $p_m \simeq 0.497$ \\
    \end{tabular}
    \caption{Common strategy $m \pi_{n-1}$ for various values of $m$}
    %\label{fig:fbh_3h}
\end{figure}

Since $\varepsilon$ takes values in $\{0,1\}$, we have: 
\begin{align*}
W_{n} &\geq \int_{[0,1]^{n}} \prod_{i=1}^{n} \varepsilon\left(f_{i,m}(x_{-i})x_{i}\right) \text{d} \Lambda^{n} \\
&= \int_{[0,1]^{n}}\varepsilon(mx_{1}\dots x_{n}) \text{d} \Lambda^{n}\\
&= \int_{[0,1]^{n-1}}\left(\underbrace{\int_{[0,1]} \varepsilon(mx_{1}\dots x_{n}) \text{d} x_{1}}_{g_{m}(x_{2},\dots, x_{n})} \right) \text{d} x_{2}\dots \text{d} x_{n}
\end{align*} 

But for any $(x_{2}, \dots x_{n}) \in ]0,1]^{n-1}$, we have $g_{m}(x_{2}, \dots x_{n}) \underset{m \to \infty}{\longrightarrow} 1/2 $ by Lemma~\ref{epsilonto1/2}. Furthermore, $g_{m} \leq 1$. Hence, by the dominated convergence theorem: \begin{center}
    $\int_{[0,1]^{n}} \prod_{i=1}^{n} \varepsilon(f_{i,m}(x_{-i})x_i) \text{d} \Lambda^{n} \longrightarrow 1/2$
\end{center}

Since $W_{n} \leq 1/2$ by Lemma \ref{integral epsilon <= 1/2}, this is sufficient to conclude.
\end{proof}

This solves Levine's continuous hat game. We now see that for all $\varepsilon > 0$, it is possible to choose a strategy with a probability of success greater than $(1/2) - \varepsilon$. By enlarging the players' set of choices, we have allowed them to perfectly correlate their strategies in order to overcome the known bounds in the initial game. We also note that by giving them this choice, computing the optimal probability becomes tractable as the problem reduces from a combinatorial problem to an analytic one. \\

Despite providing a complete solution to the continuous game, the imaginary strategy 
framework does not immediately yield improved upper bounds for $V_n$. The fundamental 
obstacle is that $W_n = 1/2 > V_n$, which shows that $\widehat{\mathcal{M}}_{n-1}$ 
is too large a class to recover $V_n$ as a supremum. Obtaining upper bounds on $V_n$ 
via this framework would require either identifying a restricted subclass of imaginary 
strategies whose optimal probability equals $V_n$, or proving that no strategy in the 
mixed class $\widehat{\mathcal{S}}_1 \times \widehat{\mathcal{M}}_1$ of 
Conjecture~\ref{new conjecture} can exceed $7/20$. We discuss this direction in 
Section~5.4. \\

\subsection{Reducing the set of strategies to standard-imaginary strategies}

We have just proved that $W_{n} = 1/2$ for all $n \geq 2$. This shows that the approximation of the first problem by imaginary strategies is not tight enough. That is, $\widehat{\mathcal{M}}^{n}_{n-1}$ is too large compared to $\widehat{\mathcal{S}}^{n}_{n-1}$ to obtain useful bounds on $V_{n}$. Therefore, one should consider a set of measurable functions of smaller size for which the optimal probability is computable. \\

We return to the case where $n=2$. One possibility would be to force one of the two players to use an $\infty$-strategy while the other can freely use any imaginary strategy. Stated differently, we consider the set\footnote{Here $2^{\widehat{\mathcal{S}}_{1}}$ denotes the set of functions $2^k$ for $k \in \widehat{\mathcal{S}}_{1}$} of strategies $2^{\widehat{\mathcal{S}}_{1}} \times \widehat{\mathcal{M}}_{1}$. It is indeed a restriction of the set of strategies since:
$$\left(2^{\widehat{\mathcal{S}}_1}\right)^2 \subsetneq 2^{\widehat{\mathcal{S}}_1} \times \widehat{\mathcal{M}}_1 \subsetneq \widehat{\mathcal{M}}_1^2.$$

Following many algorithmic attempts, we conjecture that this new approximation is much better: we have found no example of such strategies yielding a probability of victory higher than $7/20$. In particular, it is natural, given the optimal strategies that we have discovered for the continuous game, to choose $x \longmapsto mx$ as imaginary strategy and $y \longmapsto \lfloor mx\rfloor$ and study the probability as $m$ goes to infinity. In fact, this strategy doesn't even achieve a probability of $7/20$. Note that $7/20$ is achievable using only $R_t$ strategies. Thus, we believe that reducing the set of strategies to standard-imaginary strategies would lead to a sharp bound. 

\begin{conjecture} We conjecture that the following holds:
    \label{new conjecture}
         $$V_2 = \underset{k,g \in  \widehat{\mathcal{S}}_1 \times \widehat{\mathcal{M}}_1 }{\sup} \int_{[0,1]^2} \varepsilon \left(2^{k(y)} x \right) \varepsilon \left(g(x) y \right)\text{d} x \text{d} y$$
\end{conjecture}

\subsection{Theorem of optimal response}

Having conjectured the previous statement, we now wish to check if it holds experimentally. More precisely, we would like to know how to maximize - with as few computations as possible - the following quantity 

$$P_{f,g} := \int_{[0,1]^2} \varepsilon(f(y)x) \varepsilon(g(x)y) \text{d} x\text{d} y$$

for some $g \in \widehat{\mathcal{M}}_{1}$ and where $f \in \widehat{\mathcal{M}}_{1}$ is fixed. Put differently, if one of the players imposes their strategy, how can we efficiently find the optimal response of their teammate ? We could also ask ourselves whether or not this maximization can be done using only $\infty$-strategies. As we will now see, this depends on the fixed strategy $f$. \\

We first consider the case of $h$-strategies $k,\ell \in \mathcal{S}_{1,h}$ where $k$ is fixed. The probability of victory is therefore:

$$ \Proba \left(U^{(k(V))}=V^{(\ell(U))}=1 \right) = \frac{1}{4^h} \sum_{i=1}^{2^h} \sum_{j=1}^{2^h} a_i^{(k(a_j))}a_j^{(\ell(a_i))} $$

Hence, it suffices to define each $\ell(a_{i})$ in order to maximize $\sum_{j=1}^{2^h} a_i^{(k(a_j))}a_j^{(\ell(a_i))}$. In other words, the best response to an $h$-strategy can be simply computed "column by column" (as defined in the visualization tool presented earlier). This idea is not as easy to generalize to imaginary strategies. We have to prove that, in order to maximize $P_{f,g}$ with $f$ fixed, it suffices to maximize, for each value of $x \in [0,1]$, the quantity $\int_0^1 \varepsilon(f(y)x) \varepsilon(uy) \text{d} y$ where $u \geq 0$ is the parameter to be optimized. \\

In what follows, we provide a proof for this result. We later make use of it to explain how to compute optimal responses to a fixed strategy. \\

\begin{theorem}
\label{optimresp}
    We fix $f \in \widehat{\mathcal{M}}_{1}$. The following equality holds:

    $$\underset{g \in \widehat{\mathcal{M}}_1}{\sup} P_{f,g} = \int_0^1  \underset{u \geq 0}{\sup} \int_0^1 \varepsilon(f(y)x) \varepsilon(uy) \mathrm{d} y \mathrm{d} x$$
\end{theorem}

\begin{proof}
For $x \in [0,1]$ and $u \geq 0$, define
\[
\phi(x,u) := \int_0^1 \varepsilon(f(y)x)\,\varepsilon(uy)\,dy, \qquad h(x) := \sup_{u \geq 0}\, \phi(x,u).
\]
The proof is in three steps. In the first one, we obtain the easier inequality $$\underset{g \in \widehat{\mathcal{M}}_1}{\sup} P_{f,g} \leq \int_0^1  \underset{u \geq 0}{\sup} \int_0^1 \varepsilon(f(y)x) \varepsilon(uy) \mathrm{d} y \mathrm{d} x$$. In the second, we prove that $h$ is measurable and finally we conclude on the desired equality.

\medskip
\noindent\textit{Step 1.}
Let $g \in \widehat{\mathcal{M}}_1$. The map $(x,y) \mapsto \varepsilon(f(y)x)\varepsilon(g(x)y)$ is measurable and bounded by $1$, so Fubini's theorem applies and gives
\[
P_{f,g} = \int_0^1 \left(\int_0^1 \varepsilon(f(y)x)\,\varepsilon(g(x)y)\,dy\right) dx
= \int_0^1 \phi(x,\,g(x))\,dx
\leq \int_0^1 h(x)\,dx.
\]
Taking the supremum over all $g \in \widehat{\mathcal{M}}_1$ yields $\sup_{g} P_{f,g} \leq \int_0^1 h(x)\,dx$.

\medskip
\noindent\textit{Step 2.}
We claim that $h(x) = \sup_{q \in \mathbb{Q}_+} \phi(x,q)$. The inequality $$\sup_{q \in\mathbb{Q}_+} \phi(x,q) \leq h(x)$$

is immediate. For the reverse, let $u_0 > 0$ (the case $u_0 = 0$ is trivial since $\phi(x,0) = 0$), and let $(u_n) \subset \mathbb{Q}_+$ be any sequence with $u_n \nearrow u_0$. \\

The set $\{y \in [0,1] : u_0 y \in \mathbb{Z}\}$ is at most countable, hence has Lebesgue measure zero. Thus for almost every $y \in [0,1]$, we have $u_0 y \notin \mathbb{Z}$, which implies $\varepsilon(u_n y) \to \varepsilon(u_0 y)$. Since $|\varepsilon| \leq 1$, the dominated convergence theorem gives $\phi(x, u_n) \to \phi(x, u_0)$, and therefore $\sup_{q \in \mathbb{Q}_+} \phi(x,q) \geq \phi(x,u_0)$. Taking the supremum over $u_0 \geq 0$ establishes the claim.

Fix an enumeration $\mathbb{Q}_+ = \{q_1, q_2, \ldots\}$. For each $n \geq 1$, the map $(x,y) \mapsto \varepsilon(f(y)x)\varepsilon(q_n y)$ is measurable and bounded, so Fubini's theorem implies that $x \mapsto \phi(x, q_n)$ is measurable. It follows that
\[
h = \sup_{n \geq 1}\, \phi(\cdot, q_n)
\]
is a countable supremum of measurable functions, hence measurable.

\medskip
\noindent\textit{Step 3.}
Fix an integer $p \geq 1$ and define, for each $n \geq 1$, the measurable set
\[
A_{p,n} := \bigl\{x \in [0,1[ \ \mid \phi(x, q_n) \geq h(x) - 1/p\bigr\}.
\]
Since $h(x) = \sup_{n \geq 1} \phi(x, q_n)$, we have $\bigcup_{n \geq 1} A_{p,n} = [0,1[$. We extract a measurable partition by setting $B_{p,1} := A_{p,1}$ and $B_{p,k} := A_{p,k} \setminus \bigcup_{j < k} A_{p,j}$ for $k \geq 2$, and we define
\[
g_p(x) := q_k \quad \text{for } x \in B_{p,k}.
\]
The function $g_p$ is measurable and takes values in $\mathbb{Q}_+ \subset \mathbb{R}_+$, hence is a non-negative measurable function on $[0,1]$, so $g_p \in \widehat{\mathcal{M}}_1$. By construction, $\phi(x, g_p(x)) \geq h(x) - 1/p$ for almost every $x \in [0,1[$, which gives
\[
P_{f,\,g_p} = \int_0^1 \phi(x,\, g_p(x))\,dx \;\geq\; \int_0^1 h(x)\,dx - \frac{1}{p}.
\]
Since this holds for every $p \geq 1$, we conclude that $\sup_{g} P_{f,g} \geq \int_0^1 h(x)\,dx$.

\medskip

Combining the upper and lower bounds yields
\[
\sup_{g \in \widehat{\mathcal{M}}_1} P_{f,g} = \int_0^1 \sup_{u \geq 0}\int_0^1 \varepsilon(f(y)x)\,\varepsilon(uy)\,dy\;dx. \qedhere
\]
\end{proof}

As can be seen from this proof, the integral framework is very well-suited for formalizing constrained optimization in the game.\\

We apply this theorem to the so-called “first black hat” strategy. In particular, we will let a player choose the first black hat strategy $f_0$ and find the optimal response of their teammate. Note that, for the black hat strategy, we have
\[
f_0(y) = 2^{-\lfloor \log_2(y) \rfloor}.
\]

\begin{proposition}
In the two-player game, the first black hat strategy is the best response to the first black hat strategy. Therefore, there exist strategies for which best responses can be found within the class of strategies $2^{\hat{\mathcal{S}}_1}$.
\end{proposition}

\begin{proof}
Let $f$ denote the first black hat strategy. We show that for every $p \geq 1$ and every $x \in I_p := \left[ \frac{1}{2^p}, \frac{1}{2^{p-1}} \right[$, we have
\[
\sup_{u \geq 0} \int_0^1 \varepsilon(f(y)x)\, \varepsilon(uy)\, \mathrm{d}y = \frac{1}{2^p} = \int_0^1 \varepsilon(f(y)x)\, \varepsilon(f(x)y)\, \mathrm{d}y.
\]

Let $x \in I_p$. For any $y \in \bigcup_{m \leq p-1} I_m = \left[ \frac{1}{2^{p-1}}, 1 \right[$, the monotonicity of $f$ yields
\[
0 \leq f(y)x < \frac{1}{2^{p-1}} \cdot f\left( \frac{1}{2^{p-1}} \right) = 1,
\]
and therefore $\varepsilon(f(y)x) = 0$. It follows that for any $u \geq 0$,
\[
\int_0^1 \varepsilon(uy)\, \varepsilon(f(y)x)\, \mathrm{d}y = \int_0^{2^{-(p-1)}} \varepsilon(uy)\, \varepsilon(f(y)x)\, \mathrm{d}y 
\leq \int_0^{2^{-(p-1)}} \varepsilon(uy)\, \mathrm{d}y 
\leq \frac{1}{2^p}.
\]

For the reverse inequality, observe that for every $y \in \left[0, \frac{1}{2^{p-1}} \right[$, we have $2^p y < 2$, which implies
\[
\varepsilon(2^p y) = 1 \iff y \in \left[ \frac{1}{2^p}, \frac{1}{2^{p-1}} \right[ = I_p.
\]

Moreover, for every $y \in I_p$, we symmetrically have $\varepsilon(f(y)x) = \varepsilon(2^p x) = \varepsilon(2^p y) = \varepsilon(f(x)y) = 1$. Hence,
\[
\int_0^1 \varepsilon(f(y)x)\, \varepsilon(f(x)y)\, \mathrm{d}y = \int_{I_p} \mathrm{d}y = \frac{1}{2^p} = \sup_{u \geq 0} \int_0^1 \varepsilon(f(y)x)\, \varepsilon(uy)\, \mathrm{d}y.
\]

That is
\[
\sup_{g \in \widehat{\mathcal{M}}_1} P_{f,g} = P_{f,f} = \frac{1}{3}.
\]
\end{proof}

It follows that for this specific strategy $f_0$, an optimal response can be obtained simply by taking $g = f_0$. In particular, the best response in this case is an $\infty$-strategy. \\

This observation supports the conjectured equality Conjecture~\ref{new conjecture}. Although the use of imaginary strategies may initially appear too coarse, it is not necessarily so when one of the players employs an $\infty$-strategy. \\

It is important to keep in mind that what has been done here for the FBH strategy can also, in theory, be done for any other fixed strategy. This shows that the analytical framework we have just presented is particularly useful for finding optimal responses to strategies in Levine's hat problem. An interesting open question using these results is whether or not it is possible to find Nash equilibria by iterating the optimal response theorem.

\section{Conclusion}

In this paper, we have studied Levine's hat problem from a new analytical 
and geometric perspective. By representing hat stacks as real numbers in 
$[0,1]$ via their binary expansion, we derived a new integral formulation 
for $V_n$ that unifies the finite and infinite cases and provides a natural 
visualization tool for strategies. \\

Building on this framework, we constructed the recursive strategy 
$\mathscr{S}_5$, which achieves the conjectured optimal winning probability 
$7/20$ for two players. Although this matches the performance of the known 
strategy $\mathscr{S}_3$, the existence of $\mathscr{S}_5$ contradicts the 
expectation from \cite{article 3} that recursive strategies of order $t > 3$ 
bring no improvement. More importantly, it yields a strictly better geometric 
convergence rate for $V_{2,h}$, and a new lower bound $U_3(p)$ for $V_2(p)$ 
that improves all previously known bounds for $p < 0.312$. The broader 
construction of $R_t$-type strategies for all odd $t \geq 3$ (Theorem~\ref{thm:Rt_exists}) 
shows that the convergence rate of $V_{2,h}$ to $7/20$ can be made 
arbitrarily close to geometric with ratio $4$. \\

Finally, we introduced and completely solved a continuous generalization of 
the problem, proving that the optimal winning probability in this setting 
equals $1/2$ for all $n \geq 2$. While this shows that imaginary strategies 
alone are too permissive to bound $V_n$ directly, the Optimal Response 
Theorem and the mixed framework $\widehat{\mathcal{S}}_1 \times 
\widehat{\mathcal{M}}_1$ suggest a tractable reformulation of the original 
problem, as expressed in conjecture~\ref{new conjecture}. This result fits into a broader pattern : every known tractable variant of Levine's hat problem relaxes either the rules of the game (biased coins, same-color winning condition) or the information structure. Our continuous generalization shows that relaxing the data structure — the binary encoding of hats — is sufficient. In this sense, $W_n = 1/2$ provides a diagnosis of where the difficulties of the problem lives. \\

We hope that this new point of view on the puzzle can open new doors and encourage research on the topic.

%%%%%%%%%%%%%%%%%%%%%%%%%%%%%%%%%%%%%%%%%%%%%%%%%
% The following optional unnumbered section is where you put personal acknowledgements,
% research grant support, and similar things.  Do not put them on the front page.
\subsection*{Acknowledgements}

We thank sincerely Lucas Gerin for his proposal of this research subject and his supervision. Furthermore, we thank Henry Swanson~\cite{article 7} for reaching out to us while this paper was an early preprint and for sharing the key ideas underlying the construction of $R_t$-type strategies. We are also grateful to Kuo-Han Ku from the National Center for Theoretical Sciences for contacting us and presenting his work on the problem during the same period.

%BIBLIOGRAPHY
% You do not have to use the same format for your references, but 
%    include everything in this file.
% If you use BibTeX to create a bibliography, copy the .bbl file into here.
% We recommend you use \doi{...} and \arxiv{...} like the examples below,
% as they give a short display form with an active link to the full url.

\end{document}